\documentclass[journal]{IEEEtran}
\pdfoutput=1
\usepackage{ifpdf}
\usepackage[ruled]{algorithm2e}
\usepackage{color}
\usepackage[caption=false]{subfig}
\usepackage{url}
\usepackage{graphicx}
\usepackage{epstopdf}
\usepackage{float}
\usepackage{tabularx}
\usepackage{epstopdf}
\usepackage{comment}
\usepackage{amsmath,amsthm,amsfonts,amssymb}
\usepackage[noadjust]{cite}

\theoremstyle{definition}
\newtheorem{theorem}{\textbf{Theorem}}

\newtheorem{remark}{\textbf{Remark}}
\newtheorem{lemma}{\textbf{Lemma}}
\newtheorem{corollary}{\textbf{Corollary}}
\newtheorem{proposition}{Proposition}

\usepackage{tikz}
\usetikzlibrary{automata,positioning}

\newcommand{\bd}{\mathbf}
\newcommand{\ud}[1]{_\mathrm{#1}}
\newcommand{\up}[1]{^\mathrm{#1}}
\newcommand\Tstrut{\rule{0pt}{2.6ex}}         
\newcommand\Bstrut{\rule[-0.9ex]{0pt}{0pt}}   

\newcommand{\revise}[1]{#1}
\newcommand{\highlight}[1]{\colorbox{pink}{#1}}

\usepackage{lipsum}

\title{Optimal Battery Control Under Cycle Aging Mechanisms in Pay for Performance Settings}
\author{Yuanyuan Shi, Bolun Xu, Yushi Tan, Daniel Kirschen, Baosen Zhang
	\thanks{The authors are with the Department of Electrical Engineering,
		University of Washington,
		Seattle, Washington 98125,
		(e-mail:\{yyshi, xubolun, 	ystan, kirschen, zhangbao\}@uw.edu) The authors are partially supported by the University of Washington Clean Energy Institute.}
}

\IEEEoverridecommandlockouts

\begin{document}
\maketitle
\begin{abstract}
\label{sec:abstract}
We study the optimal control of battery energy storage under a general ``pay-for-performance'' setup such as providing frequency regulation and renewable integration. \revise{In these settings, batteries need to carefully balance the trade-off between following the instruction signals and their degradation costs in real-time.} Existing battery control strategies either do not consider the uncertainty of future signals, or cannot accurately account for battery cycle aging mechanism during operation. In this work, we take a different approach to the optimal battery control problem. Instead of attacking the complexity of battery degradation function or the lack of future information one at a time, we address these two challenges together in a joint fashion. In particular, we present an electrochemically accurate and trackable battery degradation model called the rainflow cycle-based model. We prove the degradation cost is convex. Then we propose an online control policy with a simple threshold structure and show it achieve near-optimal performance with respect to an offline controller that has complete future information. \revise{We explicitly characterize the optimality gap and show it is independent to the duration of operation.} Simulation results with both synthetic and real regulation traces are conducted to illustrate the theoretical results.
\end{abstract}

\section{Introduction}
\label{sec:intro}
A confluence of industry drivers -- \revise{including increased deployment of renewable generations}, the
high capital cost of managing grid peak demands, and large investments in grid infrastructure for
reliability -- has created keen interest in building and employing more energy storage systems~\cite{rastler2010electricity}. Plenty of energy storage technologies have been developed to serve different applications, such as pumped hydro-power, compressed air energy storage, batteries, flywheels and many more~\cite{dunn2011electrical}. \revise{Among these different technologies, \emph{battery energy storage (BES)} (e.g., lithium-ion batteries) features quick response time, high round-trip efficiency, pollution-free operation, and flexible power/energy ratings~\cite{xu8011541}.} These characteristics make it an ideal choice for a wide range of power system applications, including integration of renewable resources~\cite{bitar2011role}, grid frequency regulation~\cite{shi2016leveraging} and behind-the-meter load management of commercial and residential users~\cite{li2011optimal}.  For example, in 2015, there are 153.5MW newly installed battery energy storage devices in the US~\cite{gridscalebat2016}, which is roughly four times the amount of BES installment in 2014. \revise{It's worth mentioning that over 80\% of the installed capacity in 2015 occurred within the territory of PJM Independent System Operator (ISO), and the predominant use was frequency regulation.}

The focus of the paper is on the optimal control of battery energy storage under a general ``pay for performance'' setup: a battery is incentivized to follow certain instruction signals and is penalized when it cannot. For example, a battery participating in frequency regulation would receive a signal and is paid based on how well it follows the signal. Another important application that falls under this setup is a battery used by customers with onsite renewable generations, where the customers may need to purchase more expensive power from the grid if the battery cannot smooth out the local net demand. The common theme of the problems under the pay for performance setup is that the signal the battery should follow is inherently \emph{random} and the control decisions must be made in \emph{real-time}. Furthermore, battery storage naturally couples the decisions across time because of its finite energy and power capacities. Therefore, finding the optimal control policy for a battery is essentially a constrained online stochastic control problem~\cite{xu2017optimal}.


This online problem is challenging for two main reasons: 1) battery degradation and 2) lack of future information. \revise{A vital aspect of energy storage operation is to accurately model the operational cost of battery, which mainly comes from battery cells losing their energy capacities under repeated cycling~\cite{xu2016modeling}.} \revise{Analogous to cell phone batteries losing their capacity after several years of use, larger batteries used in the grid also lose their capacity with every charge and discharge action~(sometimes called capacity fading)~\cite{AroraEtAl1998}.} In fact, overly aggressive use of batteries can often deplete their useful capacities in a matter of months. However, battery degradation is a complex process governed by electrochemical reactions and \revise{depends on multiple environmental and utilization factors}. The second challenge of the lack of information is common to all stochastic control problems. At any given time, a decision must be made without knowing the future signals. \revise{This is further complicated by the coupling constrains introduced by battery.}

These two challenges are illustrated well in the the fast frequency regulation problem. Frequency regulation is a mechanism used by power system operators to correct the short timescale imbalance between generation and demand in the overall grid. In fast regulation (e.g., regD in PJM), a signal representing the imbalance is broadcasted every 2 or 4 seconds. Having enough energy to follow the regulation signals is critical to the function of the power system, especially as renewables increase the uncertainties in both generation and supply. Frequency regulation is also a natural application for batteries because of the fast variations and roughly zero-mean nature of the regulation signal. By participating in regulation, a battery receives a fixed payment ahead of the time. However, if it cannot follow the regulation signal, then a penalty is charged based on the mismatch. Therefore, at every time step, a battery must balance its degradation from following the regulation signals with the penalty of not doing so, while not knowing the future value of the signal.

In the past, many studies have attacked the battery control problem by focusing on one of the challenges. On one hand, by assuming the degradation of batteries is a quadratic function of the charge/discharge powers, we recover a type of constrained stochastic quadratic regulation problem where the key challenge is the lack of future information~\cite{chemali2015minimizing,BouchardEtAl2010}. On the other hand, one can focus on the degradation of the batteries, by employing accurate electrochemical models while assuming full knowledge of the future~\cite{northrop2011coordinate,he2016optimal}. Both directions have led to significant advances by still remain unsatisfactory. Even by assuming the signal that a battery faces is Gaussian, a constrained linear quadratic Gaussian problem is still extremely challenging to solve and provide any theoretical performance guarantees. Similarly, solving the optimization problem with accurate electrochemical models is by no means trivial even under full knowledge, and it is usually difficult to adapt the solutions to an online form. Given these difficulties, batteries still only serve as emergency backup, or used actively in grid services when they are owned by the utilities and are subsidized under renewable portfolio incentives.

In this paper, we take a different approach to the battery control problem. Instead of attacking the complexity of the degradation function or the lack of future information one at a time, we address these two challenges together in a joint fashion. Surprisingly, we provide a provably near optimal online algorithm for battery control. In particular, we show that under a form of so-called cycle based degradations, there is an \emph{online strategy that is within a constant additive gap of the optimal offline strategy under all possible future signals}. We explicitly characterize this gap and relate it to the set of possible future signals.  The key insight of this result comes from a better understanding of the degradation of electrochemical batteries and how it relates to the control problem. At a high level, capacity fading of these batteries due to charging and discharging is similar to the fatigue process of materials subjected to cyclic loading~\cite{laresgoiti2015modeling}. For each cycle, the capacity fades as a function of the depth of that cycle. In past approaches, these cycles were studied in the time domain, leading to complex optimization problems. In contrast, we look at the problem in the \emph{cycle-domain}, where the problem naturally decouples according to each cycle of the charge/discharge profile. This approach has a loose analogy with time/frequency duality, where some problems are much simpler in the frequency domain than in the time domain. Altogether, our work makes three contributions to the current state-of-art in battery control:
\begin{enumerate}
	\item  We prove the convexity of the rainflow cycle counting algorithm, which enables this electrochemically accurate model to be used in various battery optimization problems and guarantees the solution quality.
	\item We provide a subgradient algorithm to solve the optimization problem efficiently and optimally for offline battery planning and dispatch.
	\item We offer an online battery control policy with a simple threshold structure, and achieve near-optimal performance with respect to a offline controller that has complete future information.
\end{enumerate}

The online control policy proposed in this paper takes a simple threshold structure which limits a battery's state of charge~(SoC).  \revise{It reacts to new battery instructions without having to solve new optimization problems, leading to better computational performances than algorithms based on model predictive control and dynamic programming.} Compared with traditional threshold control strategies, such as pre-fixed SoC bounds~\cite{oudalov2007optimizing}, or proportional integral (PI) controller~\cite{borsche2013power}, our policy incorporates the application market prices and battery aging model into the SoC threshold calculations, which improves the model accuracy and making it applicable to most electrochemical battery cells. These considerations allows us to derive performance guarantees in form of a bounded constant gap to the full information optimal solution. \revise{This battery control policy can be applied to any power system application that faces stochastic signals and has constant prices over a specific period}, such as frequency regulation and behind-the-meter peak shaving.

The rest of the paper is organized as follows. Section \ref{sec:prior} covers the background and prior works on the battery control problem. Section \ref{sec:model} describes the proposed rainflow cycle-based degradation model. Section \ref{sec:convexity} sketches the convexity proof and the subgradient algorithm for solving the offline problem. Section~\ref{sec:policy} describes the proposed online control strategy and the optimality proof. We provide a case study in Section \ref{sec:simulation} using real data from PJM frequency regulation market, and demonstrate the effectiveness of the proposed control algorithm in maximizing profits as well as extending battery lifetime. Finally, Section \ref{sec:con} concludes the paper and outlines directions for future work.

A preliminary version of this paper has appeared in~\cite{ShiCDC2017}. The current paper expands significantly on~\cite{ShiCDC2017} and in particular, Theorem~\ref{theo1} on the convexity of the cycle depth degradation function and the simulations are new.

\section{Background and Prior Works} \label{sec:prior}
The operation of battery energy storage has received much recent research attention because of the importance of batteries to a power system with high penetration of renewables and maturing technologies~\cite{bitar2011role,qin2016online,li2011optimal,shi2016leveraging,akhil2013doe,zakeri2015electrical,xu2017optimal}. In these works, the degradation cost of the batteries are modeled in different ways.
The authors of~\cite{bitar2011role,qin2016online,xu2017optimal} assume battery has a fixed lifetime and ignore the degradation cost in optimization. This assumption works well when batteries are used sparingly, but tend to lead to overly aggressive actions for finer time resolution applications such as frequency regulation. Other energy storage control studies include degradation models either based on battery charging/discharging power~\cite{li2011optimal,shi2016leveraging} or energy throughput~\cite{akhil2013doe,ortega2014optimal,wang2017improving}. For example,~\cite{li2011optimal} assumes a convex degradation cost model based on battery charging/discharging power for households demand response, and~\cite{akhil2013doe} assign a constant price $2 \$/MWh$ based on battery energy throughput. These degradation models are convenient to be incorporated in existing optimization problems, at a cost of losing accuracy in quantifying the actual degradation cost. For example, a Lithium Nickel Manganese Cobalt Oxide (NMC) battery has \emph{ten} times more degradation when operated at near 100\% cycle depth of discharge (DoD) compared to operated at 10\% DoD for the \emph{same} amount of charged power or energy throughput~\cite{ecker2014calendar}. However, the impact of cycle depth is difficult to capture using power or energy based degradation functions.

The battery aging process is fundamentally described by a set of partial differential and algebraic equations~\cite{northrop2011coordinate,ramadesigan2012modeling}, however, they are in some sense too detailed to be used in power system applications. Even with dedicated state-of-the-art algorithms, these equations take several seconds to solve, making them too slow to be used in applications like frequency regulation where one receives a signal every 2 or 4 seconds. To mitigate these difficulties, we use a semi-empirical degradation model that combines theoretical battery aging mechanism with experimental observations. This model is motivated by viewing battery capacity fading as a material fatigue process, where a deeper charge/discharge cycle stresses battery much more than an equivalent number of shallower cycles.

\begin{figure}[!ht]
	\centering
	\includegraphics[width= 0.8 \columnwidth]{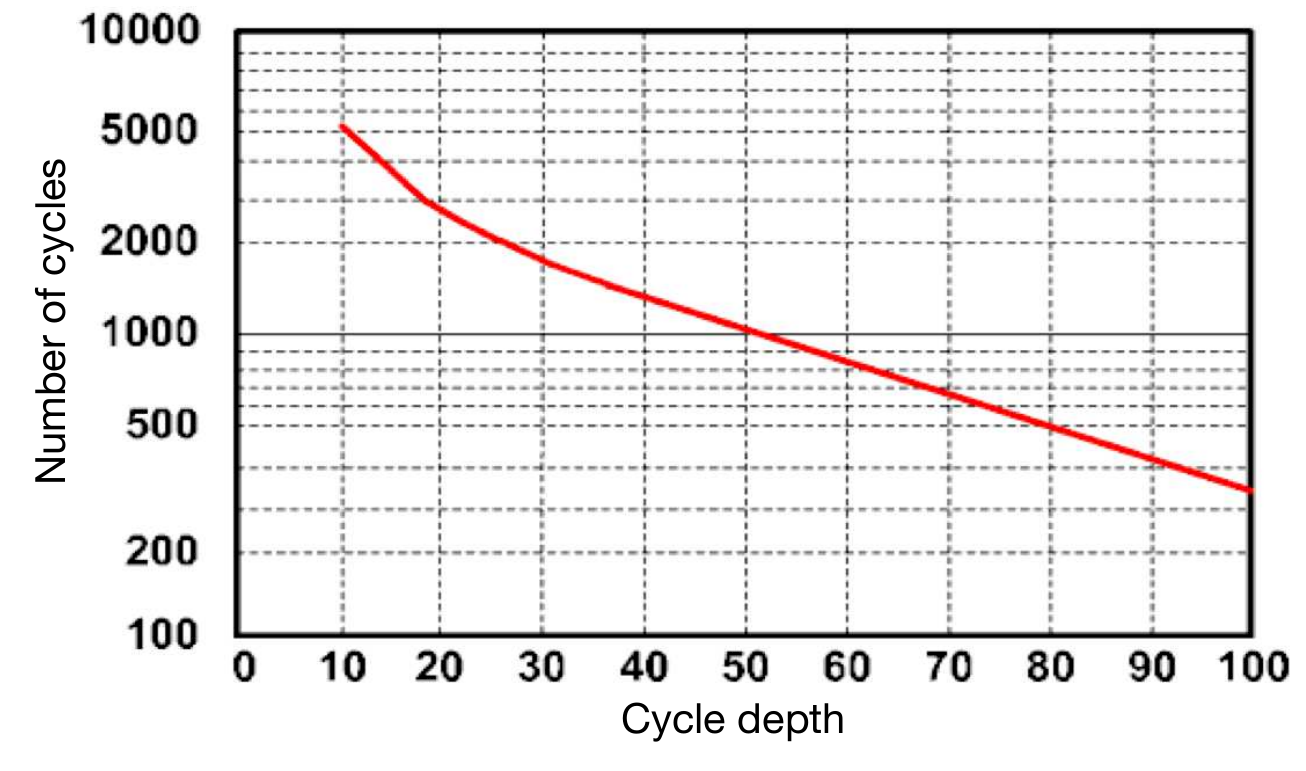}
	\caption{Battery cycle depth and operating number curve~\cite{batt_cyclecurve}. The x-axis is the cycle depth in percent, and y-axis is the number of cycles that battery could be operated under certain condition before the end of life.}
	\label{fig:depth_num}
\end{figure}

Then the relationship between cycle depth and battery degradation is defined by the cycle depth-number curve (Fig. \ref{fig:depth_num}), which are normally provided by battery manufacturers or can be estimated from field measurements.

Under this cycle aging model, each cycle causes independent stress, and the loss of battery life is the the accumulation of degradation from all cycles. A natural question is how one should count the number of cycles in a general profile, since all of them would be of heterogenous depth. Here we use the ``rainflow'' algorithm~\cite{downing1982simple}, which is the most widely adopted algorithm for cycle identification in material fatigue analysis~\cite{amzallag1994standardization} as well as for battery degradation~\cite{muenzel2015multi,dragivcevic2014capacity,he2016optimal,xu2016modeling,abdulla2016optimal}. We show that this electrochemical accurate degradation model is actually convex, which is a key step in deriving the online control algorithm and is of independent interest to many other battery applications.

The online nature of the battery control problem has perhaps received more attention from the control community. Multiple types of approaches have been developed, including model predictive control \cite{xie2012fast, qin2016online,qin2016distributed2}, stochastic and dynamic programming~\cite{bitar2011role,kim2011optimal,xu2014value,van2013optimal,melo2017online}. The authors of \cite{xie2012fast} derive a model predictive control (MPC)-based for battery energy storage and wind integration, although without any performance guarantees. Recent works~\cite{qin2016online,qin2016distributed2} do include results that bound the performance gap of online algorithms, but it is difficult to evaluate the quality of these bounds since they are either quite loose or depend on complicated optimization problems themselves that \revise{grow} with the time of operation. In addition, none of these bounds considers a cycle-based degradation problem. Our results in this paper provide an online algorithm with a constant gap to the offline optimal that is independent to the length of the operation time.

In addition to MPC type of algorithms, another widely used strategy is dynamic programming (DP). For example,~\cite{bitar2011role} and~\cite{kim2011optimal} consider using DP for storage operation with a co-located wind farm,~\cite{xu2014value} and~\cite{van2013optimal} for operating storage with end-user demands, and~\cite{melo2017online} for storage with demand response. However, for real-time control problems, the battery state space, action space and the instruction signal are all continuous. Standard DP discretization approaches tend to cause the dimension of the problem to grow exponentially. Also, implementing these algorithms requires the distributional information of the random instruction signal, which may not be readily available. In contrast, our algorithm does not require any distributional information.


\begin{remark}
  In this paper we focus on the impact of cycle-depth on the capacity lifetime of batteries. \revise{In addition to cycle-depth, numerous other factors contribute to capacity fading.} For example, the temperature of the battery has a dramatic influence in its lifetime. However, in grid applications, the temperature of the cells are normally controlled to be within a narrow band. Similarly, other factors such as extremely high C-rate and unbalanced battery cells either do not come into play for grid applications or are controlled by lower level power electronics~\cite{xu2016modeling}.
\end{remark}

\section{Model}
\label{sec:model}
In this section we describe the battery operation model, the rainflow cycle-based battery aging cost, and the pay for performance market setup. Then we state the main optimization problem on \emph{how to balance \revise{revenue} from frequency regulation and the degradation cost of battery in an online fashion.}
\subsection{Battery Operations}
We consider an operation defined over finite discrete control time steps $t\in \{1,\dotsc, T\}$, and each control time interval has a duration of $\tau$.\footnote{In practice, $\tau$ is set by the power electronic based battery management system, and is normally in the scale of milliseconds~\cite{PopEtAl2008}.} \revise{Let $x_t$ be the energy stored in the battery--the state of charge (SoC)--at the end of time $t$.} By convention, $x_t$ is a normalized quantity between $0$~(empty battery) and $1$~(full battery). At any time $t$, the battery can either charge with power $c_t$~(in units of kW) or discharge with power $d_t$~(in units of kW). Then its state of charge evolves according to the following linear difference equation~\cite{he2016optimal, wu2017optimal}:
\revise{\begin{align}
    x_{t} &= x_{t-1} + \frac{\tau \eta\ud{c}}{E} c_{t} - \frac{\tau}{\eta\ud{d}E} d_t\label{eq:bat_soc}.
\end{align}}
where the initial battery SoC is assumed to be known as $x_0$. $\eta\ud{c}$ and $\eta\ud{d}$ are the charging and discharging efficiency and $E$~(in units of kWh) is the rated energy capacity of battery. \revise{By convention, $\eta\ud{c}$ is between 0 and 1 and $\eta\ud{d}$ is larger than 1.} For ease of notations, we also define $\bd{x} = (x_t) \in \mathbb{R}^T$, $\bd{c} = (c_t) \in \mathbb{R}^T$, and $\bd{d} = (d_t) \in \mathbb{R}^T$, where $t=1,2,...,T$.

For a given battery, there are three types of operational constraints. The first is \revise{the limits of SoC}, where the stored energy in the battery is constrained to be within a pre-defined range. This constraint can arise from the health concerns since batteries like lithium-ion should not be charged completely full or discharged to be completely empty. It can also arise if batteries are used for other applications such as backup. In this paper, we assume that the SoC limits are given. The other two constraints on battery operation are the rate constraints on the charging and discharging powers, written as:
\begin{align*}
    \underline{x} \leq x_t \leq \overline{x} \mbox{, }
    0 \leq c_{t} \leq P \mbox{, and }
    0 \leq d_t \leq P,
\end{align*}
where $\underline{x}$ and $\overline{x}$ is the minimum and maximum SoC of the battery, respectively; $P$ is the battery power rating.

\revise{We consider an optimization problem where a battery is incentivized to follow an instruction signal $\mathbf{r} \in R^T$ from the system operator signifying the normalized net power imbalance in the system. We follow the convention that $r_t$ is normalized between $-1$ and $1$, where $r_t >0$ represents a shortage of power (battery is asked to provide power) and $r_t \leq 0$ represents an excess of power (battery is asked to absorb power).} Suppose the operation revenue is $R(\mathbf{c}, \mathbf{d}, \mathbf{r})$, a function of battery power output $\mathbf{c}, \mathbf{d}$ and the instruction signal. The operational cost comes from the battery degradation, denoted here by $f(\mathbf{c}, \mathbf{d})$,  a function of battery charging/discharging responses. The exact form of $f(\cdot)$, namely \emph{the rainflow cycle-based degradation function}, is introduced in the next section. The optimization objective is to maximize the net utility of the battery:
\revise{
\begin{subequations}\label{eq:battery_operation1}
	\begin{align}
	\max_{\mathbf{c}, \mathbf{d}} \; & R(\mathbf{c}, \mathbf{d}, \mathbf{r})- f(\mathbf{c}, \mathbf{d}) \label{Eq:PF_obj2} \\
	\mbox{s.t. } & x_{t} = x_{t-1} + \frac{\tau \eta\ud{c}}{E} c_{t} - \frac{\tau}{\eta\ud{d}E} d_t\,,\label{Eq:PF_C42}\\
	& \underline{x} \leq x_t \leq \overline{x}\,,\label{Eq:PF_C52}\\
	& 0\ \leq c_t \leq P \,,\label{Eq:PF_C62} \\
	& 0\ \leq d_t \leq P \,,\label{Eq:PF_C72}
	\end{align}
\end{subequations}
where \eqref{Eq:PF_C42} is the state evolution equation, \eqref{Eq:PF_C52} is the SoC constraint, \eqref{Eq:PF_C62} and \eqref{Eq:PF_C72} are the power constraints. Note here we may include a constraint that storage cannot charge and discharge at the same time~\cite{AlmassalkhiEtAl2015}, but it turns out that this condition will always be satisfied in our setting.}

Solving \eqref{eq:battery_operation1} has proven to be difficult for two reasons. The first is that most realistic cycle-based degradation functions are not well understood (e.g., they are not known to be convex), making the deterministic version of \eqref{eq:battery_operation1} nontrivial~\cite{xu2016modeling}. \revise{The second is that in real-time applications such as frequency response, the signal $\bd r$ is inherently random and difficult to forecast~\cite{HughesEtAl2017, wang2017improving}, while the state of the problem $x_t$ is constrained and coupled over time.} Therefore, even for relatively simple forms of $f$ (e.g. $f=\sum c_t^2+d_t^2$), there are no optimal or provably suboptimal online algorithms. The next section describes the rainflow cycle-based degradation model in detail, and the rest of the paper shows that rather surprisingly, this reality will lead to a simple provable optimal online algorithm.

\subsection{Cycle Counting via Rainflow}
\revise{To model the battery degradation cost $f(\mathbf{c}, \mathbf{d})$, we take the rainflow cycle-based method. This algorithm is used extensively in materials fatigue stress analysis to count cycles and quantify their depths and has also
been extensively applied to battery life assessment~\cite{muenzel2015multi,dragivcevic2014capacity,he2016optimal,xu2016modeling,abdulla2016optimal}.

The rainflow method identifies cycles from local extrema in battery SoC profile. Consider a SoC profile $\bd{x}$ with local extrema $s_1, s_2, \cdots$.  The principle of the rainflow cycle counting uses four
successive local extrema.
\begin{figure}
	\centering
	\subfloat{%
		\includegraphics[width=0.5\linewidth]{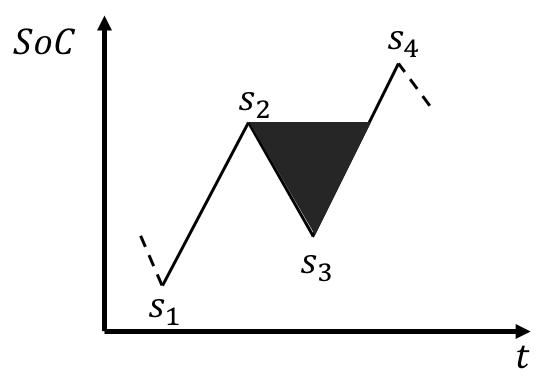}}
	\label{model: fig_1}\hfill
	\subfloat{%
		\includegraphics[width=0.5\linewidth]{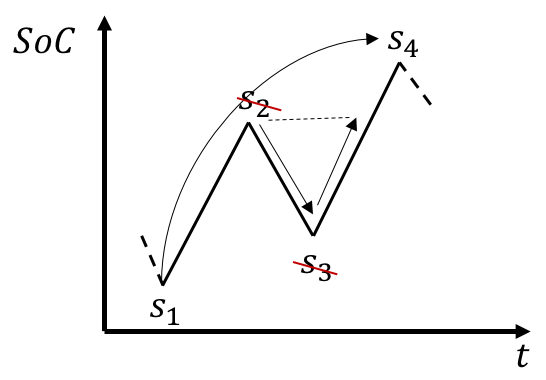}}
	\label{model: fig_2}
	\caption{Principles of rainflow cycle extraction. Given four successive point $s_1, s_2, s_3$ and $s_4$, the cycle $s_2$-$s_3$ is extracted and the points $s_2$ and $s_3$ are discarded. Then $s_1$ to $s_4$ form a charging cycle.}
	\label{fig:rainflow_procedure}
\end{figure}
This is illustrated in Figure \ref{fig:rainflow_procedure}, where $s_1, s_2, s_3$ and $s_4$ represent four successive local extrema. Three consecutive ranges are determined: $\Delta s_1 = |s_1-s_2|$, $\Delta s_2 = |s_2-s_3|$ and $\Delta s_3 = |s_3-s_4|$. If $\Delta s_2\leq \Delta s_1$ and $\Delta s_2 \leq \Delta s_3$ ($\Delta s_2$ range is less than or equal to its two adjacent ranges $\Delta s_1$ and $\Delta s_3$) then:
\begin{enumerate}
	\item a full cycle (or viewed as a charging half cycle and a discharging half cycle) represented by its extreme values $s_2$ and $s_3$ is extracted;
	\item the two points $s_2$ and $s_3$ are discarded;
	\item the two remaining parts of the sequence are connected to each other.
\end{enumerate}

If not, then the following point is considered and the same test is applied, using points $s_2, s_3, s_4$ and the new point. The procedure is repeated until the last point of the sequence is reached. After this process, the remaining points constitute what is called the \emph{residue}, in which every two consecutive points form either a charging or discharging half cycle.
\begin{figure}
	\centering
	\includegraphics[width= 0.8 \columnwidth]{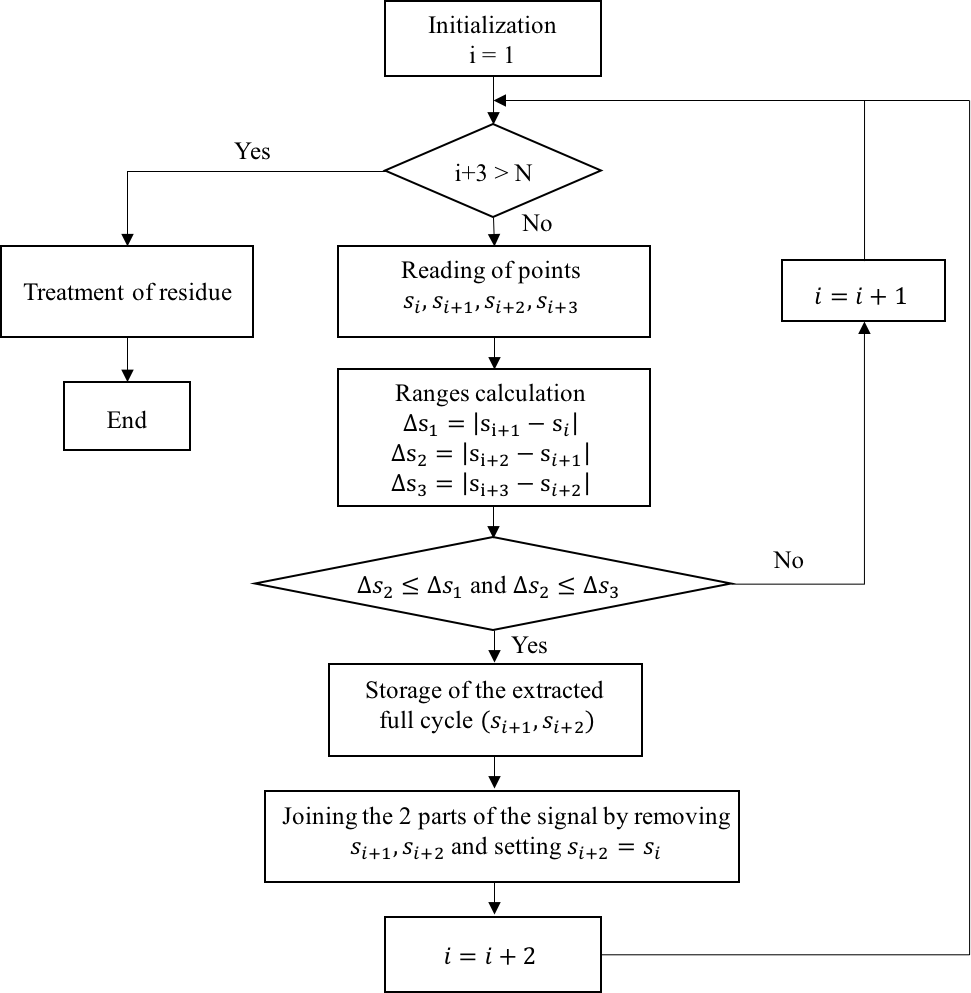}
	\caption{Flow chart of Rainflow cyle counting algorithm}
	\label{fig:flowchart}
\end{figure}

The rules of the Rainflow algorithm is summarized in flow chart Fig \ref{fig:flowchart}.  Let $\mathrm{Rainflow}$ be the functional form of the rainflow counting algorithm in Figure \ref{fig:flowchart}, where it takes an SoC profile $\bd x \in \mathbb{R}^T$ as the input, extracts all local extrema $s_1, s_2, \cdots s_N$ of $\bd x$ and outputs all the cycle depths:
\begin{align}\label{eq:rfsoc}
	(\mathbf{v}, \mathbf{w}) = \mathrm{Rainflow}(\bd x)
\end{align}
where $\mathbf{v}$ is the vector of charging half cycles and $\mathbf{w}$ is the vector of discharging half cycles.

Since cycle depths only depend on the relative differences of the turning point of the SoC profile and not on the initial SoC value, they can be calculated from $(\mathbf{c}, \mathbf{d})$
\begin{align}\label{eq:rfpower}
	(\mathbf{v}, \mathbf{w}) = \mathrm{Rainflow}\Big(\frac{\tau \eta\ud{c}}{E} \mathbf{c}- \frac{\tau}{\eta\ud{d}E} \mathbf{d}\Big).
\end{align}
}
\subsection{Battery Degradation Cost}
After counting the cycles, a cycle depth stress function $\Phi(u)$ is used to model the life loss from a single cycle of depth $u$ measured in terms of (normalized) changes in the SoC. This function indicates that if a battery cell is repetitively cycled with depth $u$, then it can operate $1/\Phi(u)$ number of cycles before reaching its end of life. \revise{The degradation cost function $\phi(\cdot)$ is normalized between 0 and 1 with respect to the total battery life. For example, if a battery can operate $100,000$ cycles before end of life at 10\% cycle depth, then a 10\% depth cycle costs $\frac{1}{100,000}$ of battery life, where $\phi(0.1) =\frac{1}{100,000}$.} In practice, this function can be estimated through empirical measurements, which is normally provided by battery manufacturers. For most electrochemical batteries,  $\Phi(u)$ is a convex function~\cite{millner2010modeling,laresgoiti2015modeling,xu2016modeling}, popularly parameterized as a power function $\alpha u^{\beta}$~\cite{xu2016modeling, laresgoiti2015modeling} or exponential functions $\alpha e^{\beta u}$~\cite{millner2010modeling}. Because cycle aging is a cumulative fatigue process~\cite{millner2010modeling,ecker2014calendar}, the total life loss $\Delta L$ is the sum of the life losses from all half cycles:
\begin{align}\label{eq:deltaL}
  \Delta L(\mathbf{v}, \mathbf{w}) = \sum_{i=1}^{|\mathbf{v}|} \frac{\Phi(v_i)}{2}+\sum_{i=1}^{|\mathbf{w}|} \frac{\Phi(w_i)}{2},
\end{align}
where $|\cdot|$ is the cardinality of a vector. Note that we assume that the stress function $\Phi$ is the same for charging and discharging, but our results hold if different functions are used.

If we substitute the rainflow algorithm as in \eqref{eq:rfpower} into \eqref{eq:deltaL}, the incremental cycle aging can therefore be written as a function of the control actions $\mathbf{c}$ and $\mathbf{d}$. To convert the loss of life to a cost, let $B$ be the battery cell replacement unit cost in \$/kWh and $E$ be the capacity of the battery in $kWh$. Then the cycle aging cost function $f(\mathbf{c}, \mathbf{d})$ is
\begin{align} \label{eqn:Jcyc}
    f(\mathbf{c}, \mathbf{d}) = \Delta L(\mathbf{c}, \mathbf{d})\cdot E \cdot B.
\end{align}

\subsection{Revenue Model}

\revise{We consider a generalized two-stage market model that captures the essence of all pay-for-performance market designs~\cite{xu2016comparison}. }  In the first stage ahead of dispatch, a payment $C$~(in units of \$) is provided to the participant. Here, we assume that this payment is known and given and focus on the second stage. The second stage occurs in real-time,
where a participant is given an instruction signal $\bd r$ and faces a penalty if it cannot follow the signal. That is, it  by pays a over-response price $\theta\in R^+$ ($\$/$MWh) for surplus injections or deficient demands during each dispatch interval, and a under-response price $\pi\in R^+$ ($\$/$MWh) for deficient injections or surplus demands. Then the total revenue is:
\begin{align}\label{eq:Regrev}
R(\mathbf{c}, \mathbf{d}, \bd r) = & C- \tau\theta\textstyle\sum_{t=1}^T|\eta_{c}c_t - \frac{d_t}{\eta_{d}} - r_t|^+ \nonumber\\
&-\tau\pi\textstyle\sum_{t=1}^T|r_t-\eta_{c}c_t+\frac{d_t}{\eta_{d}}|^+\,,
\end{align}
where $\eta_{c}c_t - \frac{d_t}{\eta_{d}}$ is the net charging power, $r_t\in [-P, P]$ is the instructed regulation dispatch set-point for the dispatch time step $t$, with the convention positive values in $r_t$ represents charging instructions.

This model captures the essence of two important applications of storage in the grid: frequency regulation\footnote{In practice, different system operators have slightly different rules for frequency response. Instead of cumbersome accounting for these rules, we focus on the generalized structure which is given in \eqref{eq:Regrev}.} and the demand shaping. \revise{In frequency regulation, $C$ is the capacity payment and $r_t$ is the regulation signal sent every 2 to 4 second by the system operator. The penalty prices $\theta$ and $\pi$ are published values. In demand shaping, a battery would enter into an agreement with a utility company to keep the demand of a customer at prescribed levels at payment $C$ and $r_t$ can be thought as the net time-varying demand of the user.} Here the penalty prices are also determined ahead of time. An important future direction is to extend our results to settings where the penalty prices are random in themselves, such as real-time arbitrage~\cite{KrishnamurthyEtAl2017,ParsonsEtAl2015}.

%

\subsection{Optimization Problem}
Summarizing the previous sections, we are attempting to solve the following optimization problem:
\begin{subequations}\label{eq:battery_operation2}
	\begin{align}
		\min_{\mathbf{c}, \mathbf{d}} \; & \tau\sum_{t=1}^T\left[\theta|\eta_{c}c_t - \frac{d_t}{\eta_{d}} - r_t|^+
		-\pi|r_t-\eta_{c}c_t+\frac{d_t}{\eta_{d}}|^+\right] \nonumber \\
		&+ \left[\sum_{i=1}^{|\mathbf{v}|} \frac{\Phi(v_i)}{2} + \sum_{i=1}^{|\mathbf{w}|} \frac{\Phi(w_i)}{2}\right] \cdot B \cdot E \label{Eq:PF_obj} \\
		\text{s.t. } & x_{t} = x_{t-1} + \frac{\tau \eta\ud{c}}{E} c_{t} - \frac{\tau}{\eta\ud{d}E} d_t \label{Eq:PF_C4}\\
		& \underline{x}  \leq x_t \leq \overline{x}\,,\label{Eq:PF_C5}\\
		& 0\ \leq c_t \leq P\,,\label{Eq:PF_C6}\\
		& 0\ \leq d_t \leq P\,.\label{Eq:PF_C7} \\
		& (\mathbf{v}, \mathbf{w}) = \mathrm{Rainflow}(\bd x). \label{Eq:PF_C2}
	\end{align}
\end{subequations}
We are interested in solving \eqref{eq:battery_operation2} in two settings: \\
\revise{
\textbf{Offline:} In the off-line setting, the entire sequence of the instruction signal $\bd r$ is given. This is important in many planning and validation problems. \\
\textbf{Online:} In the on-line setting, we solve $c_t$ and $d_t$ only based on the current and past information, $\{r_{t},r_{t-1},\dots,r_1\}$. This models the real-time decisions that batteries need to make for charging and discharging.
}

\section{Main results}
\label{sec:main}
The main contributions of this paper is to provide results that lead to optimal and computational tractable algorithms to both the offline and online solutions of \eqref{eq:battery_operation2}. The proofs of the theorems in this section are given in later sections and in the appendices. For the offline setting, we have the following theorem:
\begin{theorem}[\textbf{Convexity}] \label{theo1}
\noindent \emph{Suppose the battery cycle aging stress function $\Phi$ is convex. Then the offline version of the optimization problem in \eqref{eq:battery_operation2} is convex in the charge and discharge variables.}
\end{theorem}
This theorem settles an open question about cycle-based degradation cost functions~\cite{rychlik1996extremes,marsh2016review} and is used in the proof of the optimality of the online policy. The penalty term in the objective function \eqref{Eq:PF_obj} is clearly convex in $\bd c$ and $\bd d$, but the convexity of the term associated with the cycle stress functions is not obvious because of the nonlinear $\mbox{Rainflow}(\cdot)$ function in \eqref{Eq:PF_C2}. \revise{This result is important as it allows a range of exact approaches to convex optimization to be used on problems considering accurate battery degradation. The proof is somewhat tedious and we provide a sketch in Section~\ref{sec:theo1_sketch} and the detailed proof is found in Appendix \ref{sec:appden1}.}

\revise{Next we state the optimality result with respect to the online optimization problem. Let $J_{g}$ denote the value of \eqref{eq:battery_operation2} using the proposed online threshold control policy $g$, which will be presented later in Section \ref{sec:control_policy}. The key idea of the proposed policy is to first calculate an optimal cycle depth as a function of the degradation cost and penalty price. At each time step, the battery will follow the instruction signal until it reaches the cycle depth bound and stops following afterwards. Let $J^*$ denote the offline optima assuming all future information are known, then we have:
\begin{theorem}[ \textbf{Online optimality}] \label{theorem1}
\noindent \emph{Suppose the battery cycle aging stress function $\Phi$ is strictly convex. The proposed online control policy $g$
in Algorithm~\ref{alg:online} (given in Section~\ref{sec:control_policy}), has a constant worst-case optimality gap that is independent of the operation time duration $T$:}
 \begin{align}
	\sup_{\bd r} (J_{g}-J^*) \leq \epsilon, \text{$\forall$ $x_0$ and  $\forall$ $\{r_t\}$, $t\in\{1,\dots,T\}$.} \nonumber
 \end{align}
where $J_g$ is the cost achieved by Algorithm~\ref{alg:online}, $J^*$ is the offline optimal cost and  $\epsilon$ is a constant depending on problem parameters.
\end{theorem}}

The bound in the theorem is much tighter compared to standard bounds for online optimization problems. Normally, one would compare the averaged regret, namely $\lim_{T\rightarrow \infty} \frac{1}{T} (J_g-J^*)$ and a sublinear regret is considered to be ``good''~\cite{HazanEtAl2016,PowellEtAl2016}. Here, our result essentially shows that one can solve the online version of \eqref{eq:battery_operation2} with \emph{zero regret}, since the constant $\epsilon$ do not depend on $T$. In contrast, most existing algorithms cannot even achieve sublinear regret. Again, the key to our result is to explore the particular cyclic structure of the rainflow based cost functions. By a case study on PJM frequency regulation market in Section~\ref{sec:simulation}, we show that the proposed control algorithm could significantly improve the operational revenue up to 30\% and the battery can last as much as 4 times longer. A useful corollary of Theorem \ref{theorem1} showing when the gap is $0$:
\revise{\begin{corollary} \textbf{Zero-optimality Gap}
  \noindent\emph{If $\pi\eta\ud{d} = \theta/\eta\ud{c}$, then $J_g=J^*$. That is, there is no gap between the cost achieved by Algorithm~\ref{alg:online} and the optimal offline algorithm. }
  \label{thm:gap_0}
\end{corollary}}

For example, this corollary holds if the battery has the same charging and discharging efficiency \footnote{Again, we remind the reader that we keep the convention to write charging and discharging efficiencies differently in this paper for generality. Here, equal efficiency means $\eta_d = \frac{1}{\eta_c}$} and the penalty prices for over and under injections are the same $\theta = \pi$, then there exists an optimal online algorithm.
\revise{The proof of Theorem~\ref{theorem1} and Corollary~\ref{thm:gap_0} are given in the appendix and can be skipped if the reader wish to directly proceed to the algorithms.}

\section{Convexity and Subgradient Algorithm}
\label{sec:convexity}
In this section, we sketch the proof of Theorem \ref{theo1} to provide some intuitions and then provide the exact form of subgradient algorithm. A reader more interested in the online algorithm can directly proceed to the next section.

\subsection{Proof of Theorem \ref{theo1}}\label{sec:theo1_sketch}
Here we sketch the proof of Theorem \ref{theo1}.  Without loss of generality, we only consider the cost of charging cycles given the interchangeable and symmetric nature of charging/discharging variables. A detailed proof is given in Appendix \ref{sec:appden1}.

To prove Theorem \ref{theo1}, it suffices to show that the mapping from the SoC profile $\bd x$ to degradation cost:
	\begin{align*}
	f(\mathbf{x}) &= \left[\sum_{i=1}^{|\mathbf{v}|} \frac{\Phi(v_i)}{2} + \sum_{i=1}^{|\mathbf{w}|} \frac{\Phi(w_i)}{2}\right] \\
	 (\mathbf{v}, \mathbf{w}) &= \mathrm{Rainflow}(\bd x)
	\end{align*}
	is convex in terms of $\bd x$ given the cycle stress function $\Phi(\cdot)$ convex. That is, for any two SoC time series $\mathbf{x}, \mathbf{y} \in \mathbb{R}^{T}$,
	\begin{equation} \label{eqn:conv_f}
		f\left(\lambda \mathbf{x} + (1-\lambda) \mathbf{y}\right) \leq \lambda f( \mathbf{x}) + (1-\lambda) f(\mathbf{y}), \forall \lambda \in [0,1].
	\end{equation}
Intuitively, given two SoC series $\mathbf{x}$ and $\mathbf{y}$, if they change in different directions, the two cancel each other out so that the left hand side of \eqref{eqn:conv_f} is less than the right hand side by the convexity of $\Phi$. When $\mathbf{x}$ changes in exactly the same direction as $\mathbf{y}$ for all time steps, the equality holds. The difficulty of proving this result lies in the fact that the rainflow function is a many-to-many function that maps a sequence in $\mathbb{R}^T$ to a set of cycle depth of \emph{indeterminate length}. The proof uses induction as described in the rest of this section.

\subsubsection{Unit step decomposition}
\begin{figure}[ht]
	\centering
	\includegraphics[width=0.8 \columnwidth]{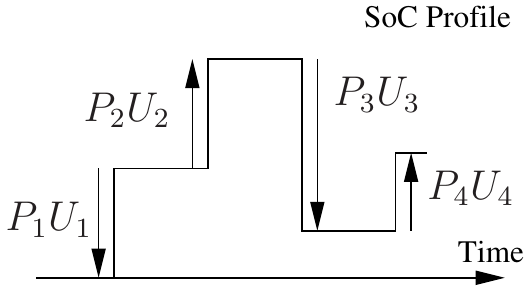}
	\caption{Decomposition of an example SoC Profile into 4 step functions.}
	\label{fig:discreted}
\end{figure}

First, we introduce the step function decomposition of SoC signal. Any SoC series $\mathbf{x}$ could be written out as a finite sum of step functions, where
\begin{align}
	\mathbf{x} = \sum_{t =1}^{T} P_t U_t \,, \label{eq:dep_s1}
\end{align}
where $U_t$ is a unit step function with a jump at time $t$ defined as:
\begin{equation*}
	U_t (\tau) =
	\begin{cases}
	1& \tau \geq t\\
	0 & \text{otherwise.}
	\end{cases}
\end{equation*}
 Fig. \ref{fig:discreted} gives an example of step function decomposition of $\bd x$.
We use this decomposition to write out $\mathbf{x}$, $\mathbf{y}$ and $\lambda \mathbf{x} + (1-\lambda) \mathbf{y}$ as finite sum of step functions, where
\begin{align}
	& \mathbf{x} = \sum_{t =1}^{T} P_t U_t\,, \mathbf{y} = \sum_{t =1}^{T} Q_t U_t \,,\label{eq:de2}\\
	& \lambda \mathbf{x} + (1-\lambda) \mathbf{y}  = \sum_{t=1}^{T} Z_t U_t\,.  \label{eq:de3}
\end{align}
Note for $\bd x$ and $\bd y$ of different length, we can take $T$ to be the maximum length since 0 can be appended to the shorter profile.

We use induction to prove Theorem \ref{theo1} on the number of non-zero step changes in $\bd y$. The base case is given in the next subsection, where $y$ has a single step change.

\subsubsection{Initial case} \label{sec:conv_initial}
\begin{figure}[ht]
		\centering
		\includegraphics[width=0.8 \columnwidth]{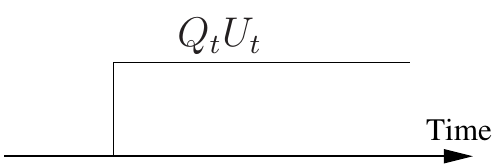}
		\caption{Base case of the induction, where one of the profiles consists of a single step.}
		\label{fig_onechange}
\end{figure}

We first show that $f(\mathbf{x})$ is convex when a profile has only one non-zero step change as shown in Fig. \ref{fig_onechange}:
\begin{lemma} \label{lemma:single_step_proof}
Under the conditions in Theorem \ref{theo1}, the rainflow cycle-based cost function $f$ satisfies
	\begin{equation*}
	f\left(\lambda \mathbf{x} + (1-\lambda)Q_t U_t \right) \leq \lambda f(\mathbf{x}) + (1-\lambda) f(Q_t U_t )\,, \forall \; \lambda \in [0,1]\,,
	\end{equation*}
where $\mathbf{x} \in \mathbb{R}^T$, and $Q_t U_t$ is a step function with a jump happens at time $t$ with amplitude $Q_t$.
\end{lemma}
The proof of this initial case requires analyzing the impact on all cycle depths from the single step and is given in Appendix \ref{sec:appden1}.

\subsubsection{Induction Steps}
Assuming Theorem 1 is true if one of the two profiles $\mathbf{x}$ or $\mathbf{y}$ has a single non-zero step. Now,  \revise{assuming} $f$ is convex up to the sum of $K$ step changes (arranged by time index):
	$$f\big(\lambda \mathbf{x} + (1-\lambda) \mathbf{y}\big) \leq \lambda f(\mathbf{x}) + (1-\lambda) f(\mathbf{y})\,, \lambda \in [0,1]$$
	if  $\mathbf{y}$ has $K$ non-zero step changes ($K < T$). We need to show $f$ is convex up to the sum of $K+1$ step changes~(i.e., $\bd y$ is of length $K+1$).
	\begin{figure}[ht]
		\centering
		\includegraphics[width=0.9 \columnwidth ]{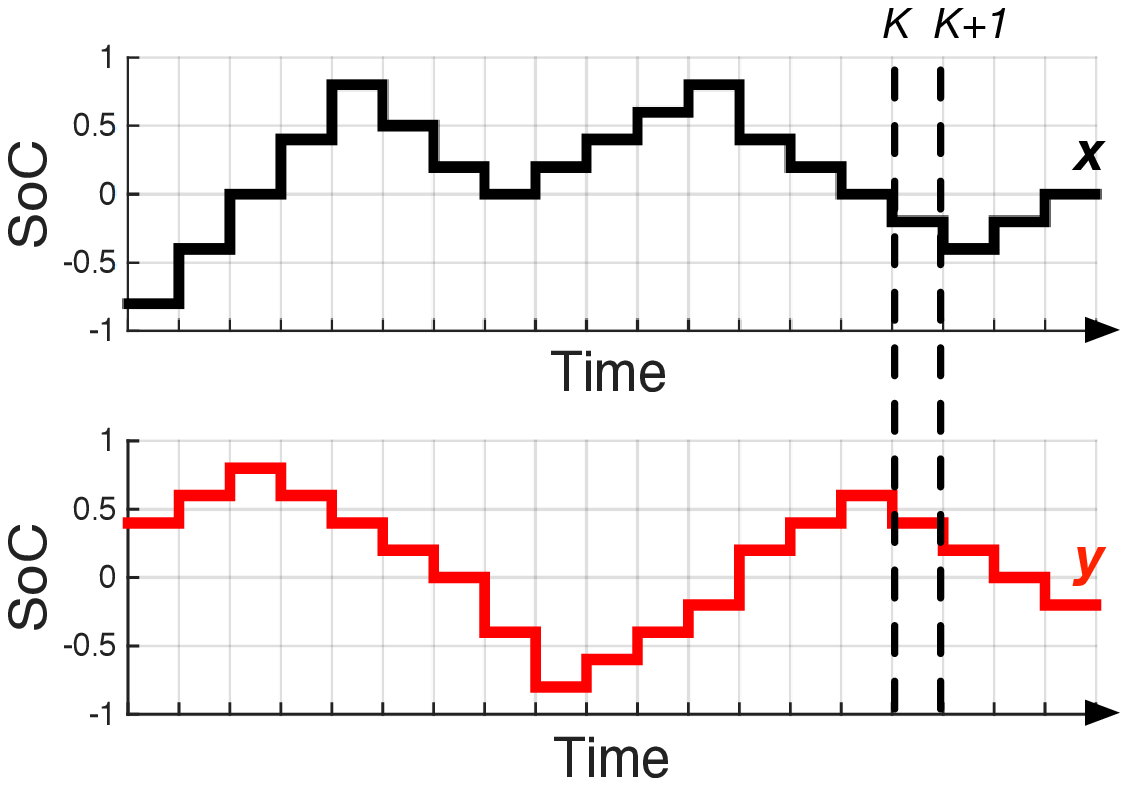}
		\caption{Induction step from $K$ to $K+1$, where we assume that $f$ is convex if one of the profiles (here the bottom profile) has only $K$ non-zero step changes and use it to show convexity of $f$ where one profile has $K+1$ non-zero steps.}
		\label{fig_induction}
	\end{figure}
The induction step proof relies on a case-by-case analysis. It contains three major conditions, where 1) $Z_{K+1}$ (the amplitude of $K+1$ step of the combined profile) and $Z_{K}$ are in the same direction 2) $Z_{K}$ and $Z_{K+1}$ are in different directions, with $|Z_{K}| \geq |Z_{K+1}|$ or 3) $Z_{K}$ and $Z_{K+1}$ are different directions, with $|Z_{K}| < |Z_{K+1}|$. Each major category may contain some further sub-cases and requires careful accounting. Showing convexity for each sub-case finishes the overall convexity proof and the detailed reasoning is given in Appendix \ref{sec:appden1}.

\subsection{Subgradient Algorithm}
\label{sec:subgradient}
\revise{The convexity of the offline problem in \eqref{eq:battery_operation2} guarantees that it can solved efficiently and optimally using gradient descend method. However, the degradation cost term $f(\bd c, \bd d)$ is not continuously differentiable (not differentiable at cycle junction points). Therefore, a subgradient method~\cite{boyd2004convex} is used in order to minimize this non-differentiable convex function. In the below section, we provide an efficient solver algorithm with the exact analytical form of subgradients. Compared with previous literatures~\cite{he2016optimal, abdulla2016optimal} using numerical solvers, the proposed algorithm with analytical subgradients form can be much faster.} With proper step size, the subgradient algorithm is guaranteed to converge to the optimal solution with a user-defined precision level~\cite{boyd2004convex}.

\revise{To begin with, we re-write the constrained optimal battery control problem in~\eqref{eq:battery_operation2} as an unconstrained optimization problem using a log-barrier function \cite{boyd2004convex}:
\begin{subequations}\label{eq:log_barrier}
	\begin{align}
		\min_{\mathbf{c}, \mathbf{d}} J(\cdot) := & \tau\sum_{t=1}^T\left[\theta|\eta_{c}c_t - \frac{d_t}{\eta_{d}} - r_t|^+
		-\pi|r_t-\eta_{c}c_t+\frac{d_t}{\eta_{d}}|^+\right] \nonumber \\
		&+ \left[\sum_{i=1}^{|\mathbf{v}|} \frac{\Phi(v_i)}{2} + \sum_{i=1}^{|\mathbf{w}|} \frac{\Phi(w_i)}{2}\right] EB  \nonumber\\
		&- \frac{1}{\lambda} \cdot \Big \{\sum_{t=1}^{T} {\log(\overline{x}-x_t)} + \sum_{t=1}^{T} {\log(x_t-\underline{x})}\nonumber\\
		&+ \sum_{t=1}^{T} {\log(P-c_t)} + \sum_{t=1}^{T} {\log(c_t)} \nonumber\\
		&+ \sum_{t=1}^{T} {\log(P-d_t)} + \sum_{t=1}^{T} {\log (d_t)} \Big\} \label{eq:log_barrier_obj} \\
		\text{s.t. } & x_t = x_0 + \sum_{k=1}^{t} \frac{\tau \eta\ud{c}}{E} c_{k} - \sum_{k=1}^{t} \frac{\tau}{\eta\ud{d}E} d_{k}\,,\label{eq:log_barrier_con1}
	\end{align}
\end{subequations}
Constraint \eqref{eq:log_barrier_con1} comes from the battery dynamics equation (8b): $x_{t} = x_{t-1} + \frac{\tau \eta\ud{c}}{E} c_{t} - \frac{\tau}{\eta\ud{d}E} d_t$.  When $\lambda \rightarrow +\infty$, the unconstrained problems \eqref{eq:log_barrier_obj} becomes equivalent to the original constrained problem.}

The major challenge of solving Eq. \eqref{eq:log_barrier_obj} lies in the second term. We need to find the mathematical relationship between charging cycle depth $v_i$ and charging power $c_t$, as well as the relationship between discharging cycle depth $w_j$ and discharging power $d_t$. Recall that the rainflow cycle counting algorithm introduced in Section~\ref{sec:model}, each time index is mapped to at least one charging half cycle or at least one discharging half cycle. Some time steps sit on the \emph{junction} of two cycles.

Let $T_{v_i}$ be all the time indexes that belong to the charging half cycle $i$ and let the time indexes belonging to the discharge half cycle $j$ be set $T_{W_j}$. Then
\begin{align}
T_{v_1}\cup\dotsc\cup T_{v_{|\bd v|}}\cup T_{w_1}\cup\dotsc\cup T_{w_{|\bd w|}} &= \{1,\dots,T\}\,,\\
T_{v_i}\cap T_{w_j} & = \emptyset\,,\, \forall i,j \label{eq:chargingordis}\,.
\end{align}

Eq. \eqref{eq:chargingordis} shows there is no overlapping between a charging and a discharging cycle. That is, each half-cycle is either charging or discharging. The cycle depth therefore equals to the sum of battery charging or discharging within the cycle time frame,
\begin{align}
v_i = \sum_{t\in T_{v_i}} \frac{\tau \eta_{c}}{E}c_t\,, \label{Eq:rf_ch}\\
w_j = \sum_{t\in TW_j} \frac{\tau}{\eta_{d}E} d_t. \label{Eq:rf_dc}
\end{align}

The rainflow cycle cost $f(\bd x)$ is not continuously differentiable. At each cycle junction point, it has more than one subgradient. We use $\partial f(\bd x)|_{c_t}$ to denote a subgradient at $c_t$. Since the SoC profile $\bd x$ is a function of $\bd c$, by the chain rule, we have
\begin{equation}
\partial f(\bd x)|_{c_t} = \Phi^{'}(v_i) \frac{B \tau \eta_{c}}{2}\,, t\in T_{v_i}\,,
\label{Eq:pardf_ch}
\end{equation}
where $v_i$ is the depth of cycle that $c_t$ belongs to. Note, at junction points, $c_t$ belongs to two cycles so that the subgradient is not unique. We can set $v_i$ to any value between $v_{i1}$ and $v_{i2}$, where $v_{i1}$ and $v_{i2}$ are the depths of two junction cycles $c_t$ belongs to.

Similarly for discharging cycle, a subgradient at $d_t$ is，
\begin{align}
\partial f(\bd x)|_{d_t} = \Phi^{'}(w_j) \frac{B \tau} {2\eta_{d}}\,, t\in TW_j \label{Eq:pardf_dc}
\end{align}
where $w_j$ is the depth of the cycle that $d_t$ belongs to. At the junction point, $w_j$ could be set to any value between $w_{j1}$ and $w_{j2}$, which are the two junction cycles $d_t$ belongs to.

Therefore, we write the subgradient of $J(\cdot)$ with respect to $c_t$ and $d_t$ as $\partial J|_{c_t}$ and $\partial J|_{d_t}$, where
\begin{small}
\begin{align}
\partial J|_{c_t} &= -\frac{\partial R(\bd c, \bd d, \bd r)}{\partial c_t} +
\Phi^{'}(v_i) \frac{B \tau \eta_{c}}{2} - \frac{1}{\lambda} \Big\{\sum_{k=t}^{T} \frac{1}{x(k)-\overline{x}}(\frac{\tau \eta_{c}}{E}) \nonumber\\
&+ \sum_{k=t}^{T} \frac{1}{x(k)-\underline{x}}(\frac{\tau \eta_{c}}{E}) + \frac{1}{c_t-P} + \frac{1}{c_t} \Big\}, t \in T_{v_i} \\
\partial J|_{d_t} &= -\frac{\partial R(\bd c, \bd d, \bd r)}{\partial d_t}  + \Phi^{'}(w_j) \frac{B \tau} {2\eta_{d}} -\frac{1}{\lambda} \Big\{\!-\sum_{k=t}^{T} \frac{1}{x(k)\!-\!\overline{x}}(\frac{\tau}{\eta_{d} E})\!\nonumber\\
& -\!\sum_{k=t}^{T} \frac{1}{x(k )\!-\!\underline{x}}(\frac{\tau}{\eta_{d} E}) + \frac{1}{d_t-P} + \frac{1}{d_t} \Big\}, t \in TW_j
\end{align}
\end{small}

The update rules for $c_t$ and $d_t$ at the $k$th iteration are,
\begin{equation*}
c_{(k)}(t) = c_{(k-1)}(t) - \alpha_k \cdot \partial J|_{c_{(k-1)}(t)}\,,
\label{eq:charging_update}
\end{equation*}
\begin{equation*}
d_{(k)}(t) = d_{(k-1)}(t) - \alpha_k \cdot \partial J|_{d_{(k-1)}(t)}\,,
\label{eq:discharging_update}
\end{equation*}
where $\alpha_k$ is the step length at $k$th iteration. Since the subgradient method is not a decent method \cite{boyd2004convex}, it is common to keep track of the best point found so far, i.e., the one with smallest function value. At each step, we set
\begin{equation*}
J_{(k)}^{best} = \min \big\{J_{(k-1)}^{best}, J(\mathbf{c}_{(k)}, \mathbf{d}_{(k)})\big\}\,,
\end{equation*}
Since the $J(\cdot)$ is convex, choosing an appropriate step size guarantees convergence.

\section{Online Policy}
\label{sec:policy}
In this section, we describe the proposed online battery control policy which balances the cost of deviating from the instruction signal and the cycle aging cost of batteries while satisfying operation constraints. This policy takes a \emph{threshold form} and achieves an optimality gap that is \emph{independent of the total number of time steps}. Therefore in term of regret, this policy achieves the strongest possible result: the regret do not grow with time. Note we assume the regulation capacity has already been fixed in the previous capacity settlement stage.
\subsection{Control Policy Formulation}\label{sec:control_policy}
\revise{The key part of the control policy is to calculate thresholds that bound the SoC of the battery as a function of the deviation penalty and degradation cost.} Let $\hat{u}$ denote this bound on the SoC and it is given by:
\begin{align}\label{eq:pol3}
    \hat{u} = \dot{\Phi}^{-1}\Big(\frac{\pi\eta\ud{d}+ \theta/\eta\ud{c}}{B}\Big) 
\end{align}
where $\dot{\Phi}^{-1}(\cdot)$ is the inverse function of the derivative of the cycle stress function $\Phi(\cdot)$.
\revise{Recall that $\pi$ and $\theta$ are the penalty prices of not meeting the instruction signal and $B$ is the price of replacing a cell in the battery. Since $\Phi$ is an increasing function, $\Phi^{-1}$ is also an increasing function. If the replacement cost is relatively small compared to the penalties($\frac{\pi\eta\ud{d}+ \theta/\eta\ud{c}}{B}$ is large), $\hat{u}$ would also be large, therefore allowing the battery a wider SoC range to operate in. On the other hand, if the replacement cost is large compared to the penalties ($\frac{\pi\eta\ud{d}+ \theta/\eta\ud{c}}{B}$ is small), $\hat{u}$ would be small, leading to a narrower range of SoC the battery would operate in.} 

\begin{algorithm}[!htb]
\SetAlgoLined
\KwResult{Determine battery dispatch point $c_t$, $d_t$}
 \tcp{initialization}
 set $\Phi\Big(\frac{\pi\eta\ud{d}+ \theta/\eta\ud{c}}{B}\Big) \to \hat{u}$, $x_0 \to x\up{max}_{0}$, $x_0 \to x\up{min}_{0}$\;
 \While{$t \leq T$}{
  \tcp{read $x_t$ and update controller}
  set $\max\{x\up{max}_{t-1}, x_t\} \to x\up{max}_t$, $\min\{x\up{min}_{t-1}, x_t\} \to x\up{min}_t$\;
  set $\min\{\overline{x}, x\up{min}_t + \hat{u}\} \to \overline{x}_t$\;
  set $\max\{\underline{x}, x\up{max}_t - \hat{u}\} \to \underline{x}_t$\;
  \tcp{read $r_t$ and enforce soc bound}
  \eIf{$r_t \geq 0$}{
   set $\min\Big\{\frac{E}{\tau\eta\ud{c}}(\overline{x}_t-x_t), r_t\Big\} \to c_t$, $0 \to d_t$ \;
   }{
   set $0\to c_t$, $\min\Big\{\frac{E \eta\ud{d}}{\tau}(x_t-\underline{x}_t), r_t\Big\}\to d_t$ \;
  }
  \tcp{wait until next control interval}
  set $t+1\to t$\;
 }
 \caption{Proposed Control Policy}
 \label{alg:online}
\end{algorithm}

\begin{figure}
	\centering
	\subfloat[Instruction (dotted line) vs. response (solid line)]{%
	\includegraphics[trim = 0mm 0mm 0mm 0mm, clip, width = .9\columnwidth]{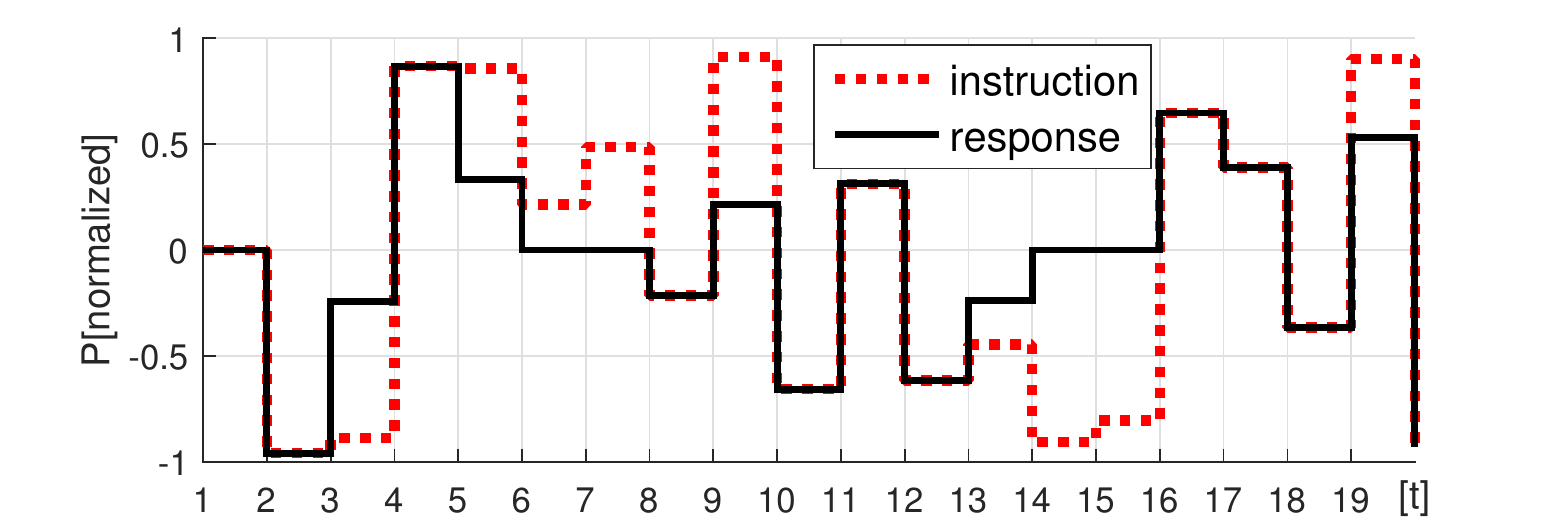}}
	\label{fig:con1}\hfill
	\subfloat[Optimal (solid) vs. profile if signal is followed perfectly (dotted)]{%
	\includegraphics[trim = 0mm 0mm 0mm 0mm, clip, width = .9\columnwidth]{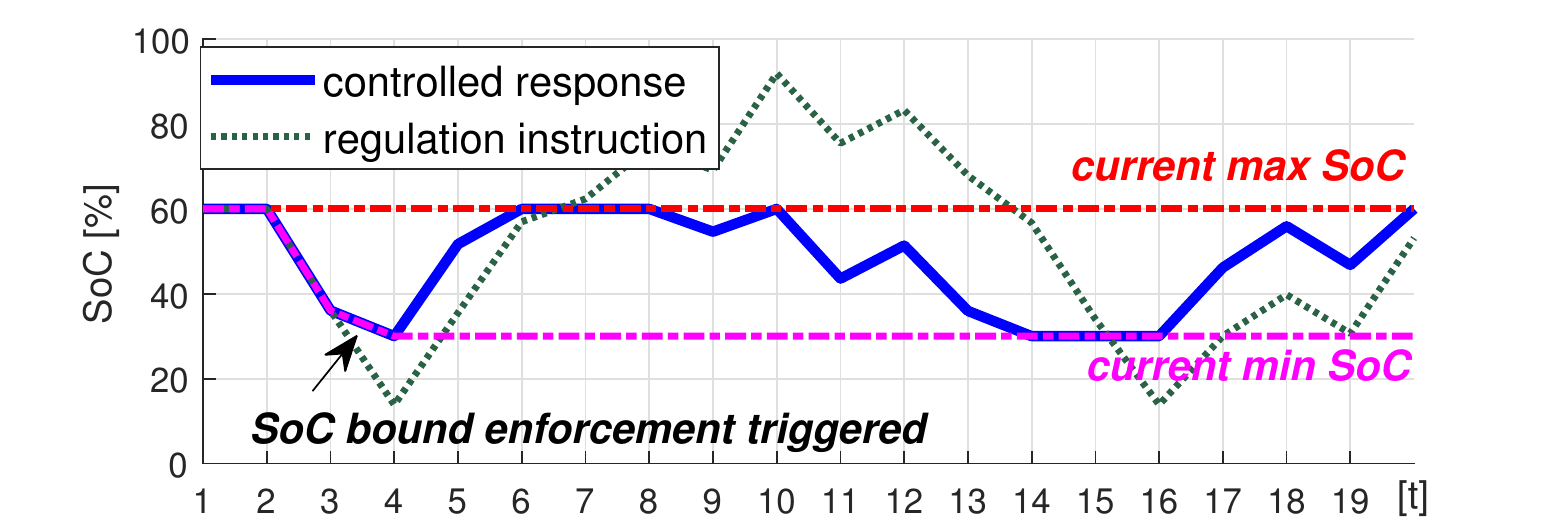}}
	\label{fig:con2}
	\caption{Illustration of the control policy in Algorithm \ref{alg:online}. The policy keeps track the current maximum and minimum SoC level. When the distance in between reaches the calculated threshold $\hat{u}$, the policy starts to constrain the response. Deeper charge and discharge cycles are avoided.}
	\label{fig:con}
\end{figure}

The proposed control policy is summarized in Algorithm~\ref{alg:online}, and Fig.~\ref{fig:con} shows a control example of the proposed policy, in which the battery follows the regulation instruction until the distance between its maximum and minimum SoC reaches $\hat{u}$.

The detailed formulation is as follows. We assume at a particular control step $t$, $x_t$ (battery state of charge) and $r_t$ (frequency regulation signal) are observed, and the proposed regulation policy has the following form: $g_t(x_t,r_t)=\begin{bmatrix} c_t & d_t \end{bmatrix}$.

The control policy employs the following strategy
\begin{align}\label{eq:pol5}
    \text{If $r_t \geq 0$, } 
		c_t &=\min\Big\{\frac{E}{\tau\eta\ud{c}}(\overline{x}_t-x_t), r_t\Big\} \\
    \text{If $r_t < 0$, } 
		d_t &=\min\Big\{\frac{E \eta\ud{d}}{\tau}(x_t-\underline{x}_t), r_t\Big\}
\end{align}
where $\overline{x}_t$ and $\underline{x}_t$ are the upper and lower storage energy level bound determined by the controller at the control interval $t$ for enforcing the SoC band $\hat{u}$
\begin{align}\label{eq:pol6}
    \overline{x}_t &= \min\{\overline{x}, x\up{min}_t + \hat{u}\}\nonumber\\
    \underline{x}_t &= \max\{\underline{x}, x\up{max}_t - \hat{u}\}\,
\end{align}
and $x\up{max}_t$, $x\up{min}_t$ is the current maximum and minimum battery storage level since the beginning of the operation, which are updated at each control step as
\begin{align}\label{eq:pol4}
    x\up{max}_t &= \max\{x\up{max}_{t-1}, x_t\}\nonumber\\
    x\up{min}_t &= \min\{x\up{min}_{t-1}, x_t\}\,.
\end{align}

\subsection{Optimality Gap to Offline Problem}\label{sec: opt_gap}
Theorem \ref{theorem1} states that the gap between the online policy in Algorithm \ref{alg:online} and an offline optimal solution is bounded by a constant.
This constant can be explicitly characterized. To do this, we define three new functions:
\begin{subequations}\label{eq:cyc_cost}
\begin{align}
    J\ud{u}(u) &=  EB\Phi(u) + E(\theta/\eta\ud{c}+\pi\eta\ud{d})u \label{eq:cyc_cost1}
    \\
    J\ud{v}(v) &= (1/2)EB\Phi(v) + (E/\eta\ud{c})\theta v \label{eq:cyc_cost2}
    \\
    J\ud{w}(w) &= (1/2)EB\Phi(w) + E\eta\ud{d}\pi w\, \label{eq:cyc_cost3}
\end{align}
\end{subequations}
where $J\ud{u}$ is the cost associated with a full cycle (made up of a charging half cycle and a discharging half cycle with equal magnitude), $J\ud{v}$ for a charge half cycle, and $J\ud{w}$ for a discharge half cycle. The detailed transforming procedure is discussed in the Appendix II-A.

If the cycle depth stress function $\Phi(\cdot)$ is strictly convex, then it is easy to see that \eqref{eq:pol3} is the unconstrained minimizer to \eqref{eq:cyc_cost1}. Similarly, the unconstrained minimizers of  \eqref{eq:cyc_cost2} and \eqref{eq:cyc_cost3} are:
\begin{align}
    \hat{v} = \dot{\Phi}^{-1}\Big(\frac{\theta/\eta\ud{c}}{B}\Big),\quad \hat{w} = \dot{\Phi}^{-1}\Big(\frac{\pi\eta\ud{d}}{B}\Big)\,.
\end{align}
The following theorem offers the analytical expression for $\epsilon$.
\begin{theorem}\label{theorem2}
If function $\Phi(\cdot)$ is strictly convex, then the worst-case optimality gap for the proposed policy $g(\cdot)$ in Theorem \ref{theorem1} is
\begin{align}
    \epsilon = \begin{cases}
    \epsilon\ud{w} & \text{if $\pi\eta\ud{d} > \theta/\eta\ud{c}$} \\
    0 & \text{if $\pi\eta\ud{d} = \theta/\eta\ud{c}$} \\
    \epsilon\ud{v} & \text{if $\pi\eta\ud{d} < \theta/\eta\ud{c}$}
    \end{cases}\,
\end{align}
where
\begin{align}
    \epsilon\ud{w} &= J\ud{w}(\hat{u})+2J\ud{v}(\hat{u})-J\ud{w}(\hat{w})-2J\ud{v}(\hat{v}) \label{eq:gap_v}\\
    \epsilon\ud{v} &= 2J\ud{w}(\hat{u})+J\ud{v}(\hat{u})-2J\ud{w}(\hat{w})-J\ud{v}(\hat{v})\label{eq:gap_u}\,.
\end{align}
\end{theorem}

\begin{figure}[ht]
	\centering
	\subfloat[Regulation instruction signal]{%
		\includegraphics[width=0.5\linewidth]{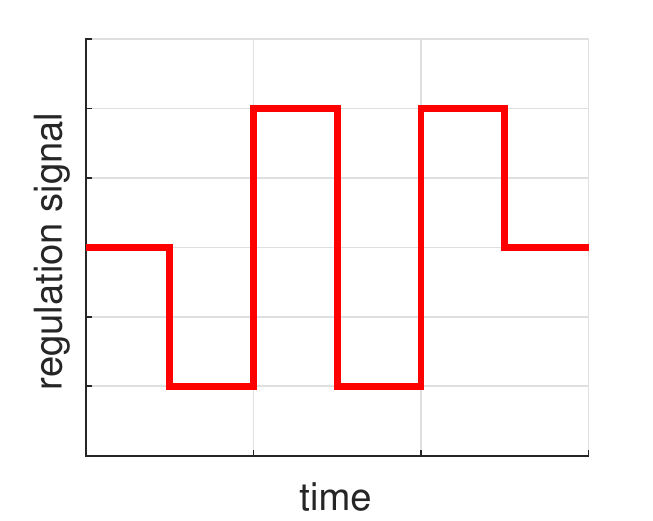}}
	\label{fig:policy_A}\hfill
	\subfloat[$\theta = \pi >0$]{%
		\includegraphics[width=0.5\linewidth]{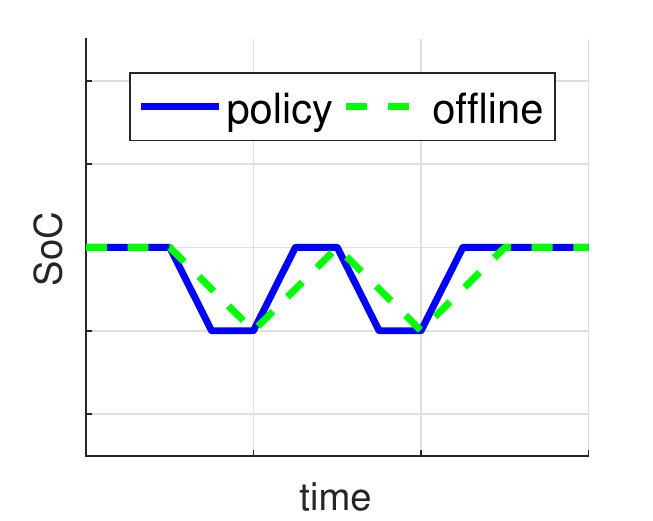}}
	\label{fig:policy_B}\\
	\subfloat[$\theta = 0, \pi>0$]{%
		\includegraphics[width=0.5\linewidth]{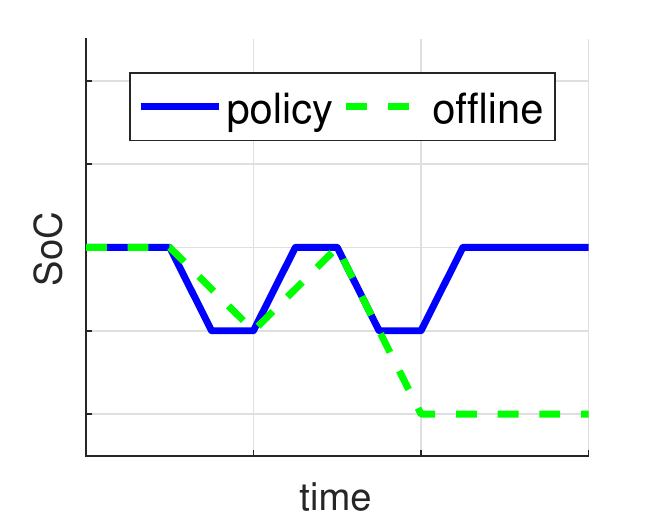}}
	\label{fig:policy_C}\hfill
	\subfloat[$\theta>0, \pi=0$]{%
		\includegraphics[width=0.5\linewidth]{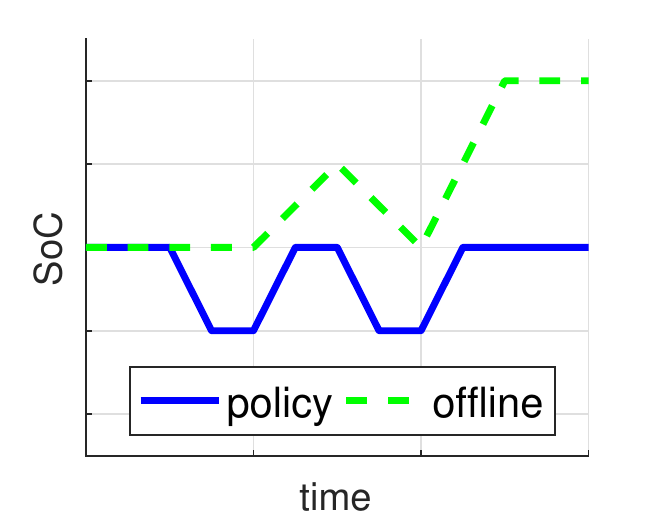}}
	\label{fig:policy_D}
	\caption{Example illustration of the policy optimality under different price settings. The value of $\theta+\pi$ is the same in all cases and the round-trip efficiency is assumed to be one, so $\hat{u}$ is the same in all cases.}
	\label{fig:policy}
\end{figure}

Note that Corollary \ref{thm:gap_0} follows from Theorem \ref{theorem2} directly. We defer the proof of the latter to Appendix \ref{sec:appden2}. The intuition is that battery operations \revise{consist mostly of} full cycles due to limited storage capacity because the battery has to be charged up before discharged, and vice versa. Enforcing $\hat{u}$--the optimal full cycle depth calculated from penalty prices and battery coefficients--ensures optimal responses in all full cycles. In cases that $\pi\eta\ud{d} = \theta/\eta\ud{c}$, $\hat{u}$ is also the optimal depth for half cycles, and the proposed policy achieves optimal control. In other cases, the optimality gap is caused by half cycles because they have different optimal depths. However, half cycles have limited occurrences in a battery operation because they are incomplete cycles~\cite{amzallag1994standardization}, so that the optimality gap is bounded as stated in Theorem~\ref{theorem2}. Fig.~\ref{fig:policy} shows some examples of the policy optimality when responding to the regulation instruction (Fig~\ref{fig:policy_A}) under different price settings. The proposed policy has the same control action in all three price settings because of the same $\hat{u}$. The policy achieves optimal control in Fig~\ref{fig:policy_B} because $\hat{u}$ is the optimal depth for all cycles. In Fig~\ref{fig:policy_C} and Fig~\ref{fig:policy_D}, half cycles have different optimal depths and the policy is only near-optimal. However, the offline result also selectively responses to instructions with a zero penalty price (charge instructions in Fig~\ref{fig:policy_C}, discharge instructions in Fig~\ref{fig:policy_D}), because it returns the battery to a shallower cycle depth with smaller marginal cost.

\newpage
\section{Simulation Results}
\label{sec:simulation}
\begin{table*}[!htb]
	\centering
	\vspace{-8mm}
	\caption{Simulation with Random Generated Regulation Signals.}
	\begin{tabular}{l c c c c c c c c c c c}
		\hline
		\hline
		& $\theta$ & $\pi$ & $\eta$ & $T$ & $\hat{u}$ & \multicolumn{2}{c}{Average objective value [\$]} & Theoretical worst-case & Maximum optimality gap  
		\Tstrut\\
		Case & [\$/MWh] & [\$/MWh] & [\%] & & [\%] &Offline & Proposed controller & optimality gap $\epsilon$ [\$] &  among 100 simulations [\$]
		\Bstrut\\
		\hline
		1 & 50 & 50 & 100 & 100 & 11.1 & 117.4 & 117.4 & \highlight{0.00} & \ \highlight{0.00} \Tstrut\\
		2 & 100 & 100 & 100 & 100 & 21.9 & 168.7 & 168.7 &  \highlight{0.00} & \highlight{0.00} \\
		3 & 200 & 200 & 100 & 100 & 42.8 & 219.4 & 219.4 &  \highlight{0.00} & \highlight{0.00} \\
		4 & 50 & 50 & 85 & 100 & 11.2 & 117.2 & 117.3 & \highlight{0.06} & \highlight{0.06}  \\
		5 & 80 & 20 & 85 & 100 & 11.7 & 108.0 & 110.7 & \highlight{3.83}  & \highlight{3.83}  \\
		6 & 20 & 80 & 85 & 100 & 10.6 & 122.4 & 123.8 & \highlight{2.19}  & \highlight{2.19}\\
		7 & 50 & 50 & 85 & 200 & 11.2 & 235.6 & 235.7 & \highlight{0.06}  & \highlight{0.06}  \\
		8 & 80 & 20 & 85 & 200 & 11.7 & 219.5 & 222.2 & \highlight{3.83}  & \highlight{3.83} \\
		9 & 20 & 80 & 85 & 200 & 10.6 & 247.6 & 248.9 & \highlight{2.19}  &  \highlight{2.19}\Bstrut\\
		\hline
		\hline
	\end{tabular}
	\label{tab:sim}
\end{table*}
\revise{
	In this section, we present several simulation examples from different aspects. In section \ref{sec:result1}, we compare the  performance of Rainflow cycle-based model against the benchmark energy throughput degradation model in offline battery operation optimization; and verified its efficiency in increasing the BES operational profit and extending battery lifetime.  In section \ref{sec:result2}, we validate the optimality of the proposed online control policy using massive random generated regulation traces. In Section VI-C, we compare the proposed online BES control policy against two state-of-the methods using realistic regulation signals from PJM Interconnection~\cite{pjm}: a greedy controller~\cite{bitar2011role} and a MPC controller assuming perfect future information is known~\cite{wang2017improving}.
}

\begin{table}[H]
	\centering
	\caption{Key Properties of the Simulated Battery}
	\begin{tabular}{p{2.8cm} p{5cm}}
		\hline
		\hline
		Attributes &   Value \Tstrut\Bstrut\\
		\hline
		Peak power  &1MW \Tstrut\Bstrut\\
		Capacity &  0.25MWh \Tstrut\Bstrut\\
		Cell type &  Lithium-ion \Tstrut\Bstrut\\
		Efficiency &  95\% for both charging and discharging \Tstrut\Bstrut\\
		Lifetime & 3,000 cycles at 80\% cycle depth \Tstrut\Bstrut\\
		Cell price & 300~\$/kWh \Tstrut\Bstrut\\
		\hline
		\hline
	\end{tabular}
	\label{tab:battery}
\end{table}
\revise{
\subsection{Comparison between the Rainflow Cycle-based Model and Benchmark Energy Throughput Model}\label{sec:result1}
To demonstrate the efficiency of the rainflow cycle-based degradation model in maximizing the BES operation revenue and extending BES lifetime, we compare with the benchmark linear energy throughput model. Linear energy throughput model is one of the most widely adopted degradation model in previous BES optimization literatures~\cite{akhil2013doe,ortega2014optimal,wang2017improving}. Consider a lithium-ion battery with a polynomial cycle depth stress function concluded from lab tests~\cite{laresgoiti2015modeling}, where $\Phi(u) = \mathrm{(5.24\times10^{-4})}u^{2.03}$, e.g., for a cycle with depth $u=0.8$, $\Phi(u) = 3.33 \times 10^{-4}$. Other key properties of the battery are summarized
in Table \ref{tab:battery}. For the linear energy throughput degradation model, we amortize the total battery cost across its lifetime energy throughput. Since the battery could operate for $3,000$ cycles at 80\% cycle depth, the lifetime energy throughput is $4.8 \times 10^{3}MWh$. Therefore, the linear cost coefficient $\lambda_e = \frac{300, 000}{4.8 \times 10^{3}}= 62.5\$/MWh$. 

We compare the performance of rainflow degradation model and linear energy throughput model in the offline optimal BES frequency regulation problem. Assume the battery bid for 1~MW symmetric regulation capacity, and the frequency regulation capacity payment is 50\$/MWh~\cite{shi2016leveraging}. The main trade-off for battery in the optimal frequency regulation problem is between the mismatch penalty and cell degradation cost. Depending on the mismatch penalty price, the optimized battery response might be quite different. We consider the following two cases: 1) when the mismatch penalty $\theta=\pi = 100$, which are higher than the linearized battery degradation cost $\lambda_e$; and 2) when the mismatch penalty $\theta=\pi = 50$, which are lower than the battery degradation cost $\lambda_e$. Fig. \ref{fig:result1}(a) and Fig. \ref{fig:result1}(b) visualize the optimized battery power outputs and SoC evolution curves in a 1 hour optimization horizon under two cases.
\begin{figure}
	\centering
	\subfloat[$\theta=\pi = 100$]{%
		\includegraphics[width=0.5\linewidth, height=5cm]{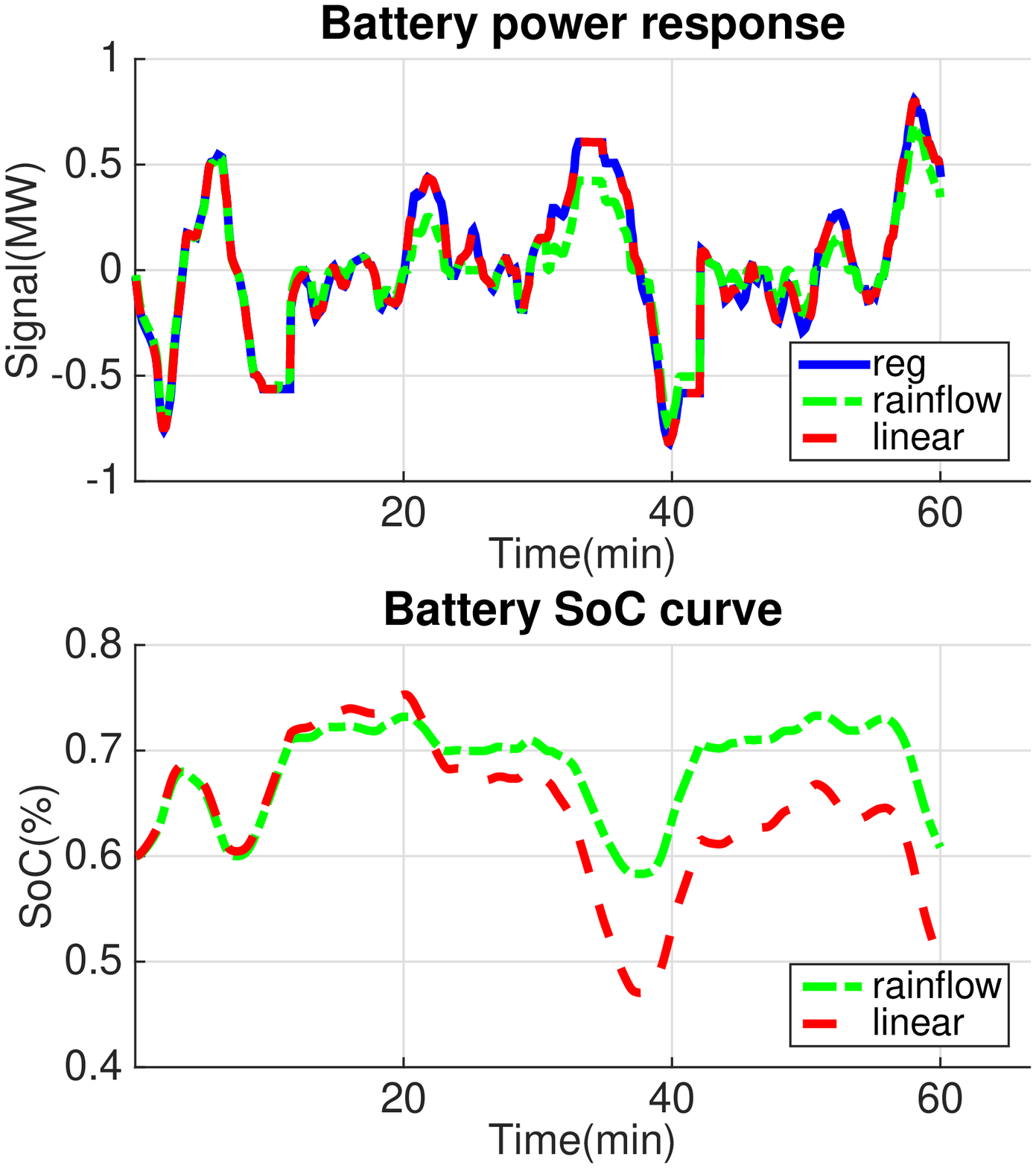}}
	\label{fig:case1}\hfill
	\subfloat[$\theta=\pi = 50$]{%
		\includegraphics[width=0.5\linewidth, height=5cm]{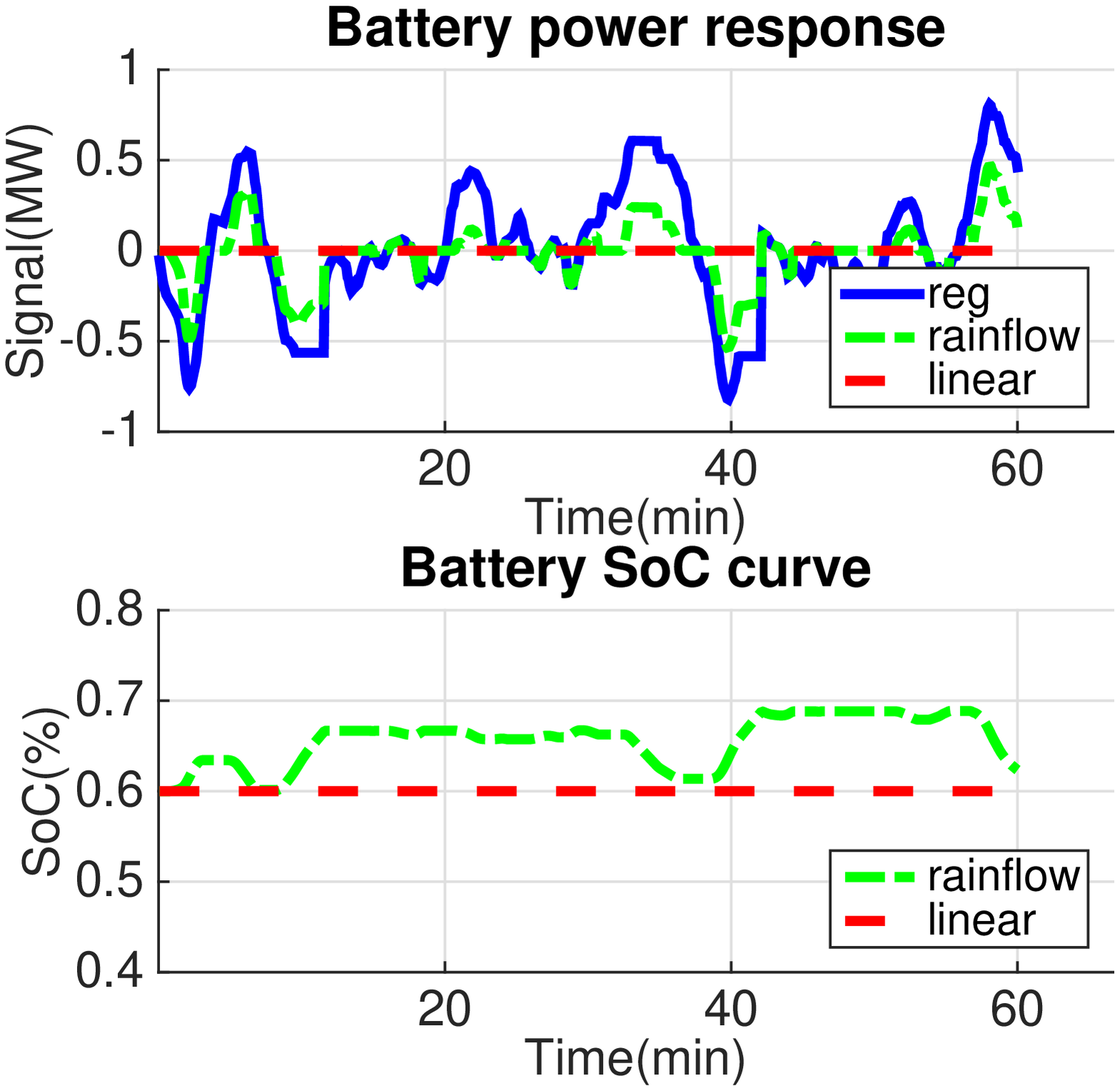}}
	\label{fig:case2}
	\caption{Comparison of the battery response under Rainflow cycle-based model and linear energy throughput model when mismatch penalty $\theta=\pi = 100$ and $\theta=\pi = 50$.}
	\label{fig:result1}
\end{figure}

In Fig. \ref{fig:case1}, when the mismatch penalty $\pi=\theta=100 > \lambda_e$, the battery completely follows the regulation instruction signal under the linear energy throughput model. However, under the rainflow cycle-based model, instead of the ``blindly following'' strategy, the battery follows the regulation signal most times and stops following the instruction signal to avoid high degradation cost of deep cycles. As we can observe from the bottom plot of Fig. \ref{fig:result1}(a), the battery SoC is restricted to a moderate range and evolves smoothly under the rainflow cycle-based cost while changes wildly under the linear model. When $\pi=\theta =50<\lambda_e$, the battery stays idle under the linear model since the linearized battery deployment cost is higher than mismatch penalty, which leads to \emph{zero} regulation service revenue. By comparison, under the rainflow model, the battery responds to the regulation signal strategically within a narrow cycle depth bound, which leads to $14.1\$$/hour revenue by solving the optimization problem.

The above comparison results show that rainflow cycle-based model can better capture the battery aging cost, therefore it provides a better BES operation schedule than previous linear energy throughput model. However, the major challenge that obstructs the adoption of rainflow cycle-based degradation model in previous BES optimization literatures is the lack of an efficient solver. Here, we compare the computational efficiency of our subgradient solver with exact gradient form in Section~\ref{sec:subgradient} with the numerical solver used in previous literatures~\cite{he2016optimal,abdulla2016optimal}. Table~\ref{tab:run_time2} shows the difference of computation time between the two solvers. It turns out that the latter does not converge for problem horizon of longer than 4 hours. All experiments conducted on a Macbook Pro with 2.5 GHz Intel Core i7, 16 GB 1600 MHz DDR3.
\begin{table}[ht]
	\centering
	\caption{Computation Time Comparison between the proposed subgradient solver with analytical subgradient form and previous numerical solvers~\cite{he2016optimal,abdulla2016optimal}}
	\begin{tabular}{l c c c c c}
		\hline
		\hline
		Time horizon (min) &  60 & 120 & 240 & 720 & 1440 \Tstrut\Bstrut\\
		\hline
		Subgradient solver time(s)  & 23.9 & 62.5 & 156.3 & 673.5 & 2522 \Tstrut\Bstrut\\
		Numerical solvers time(s)   & 264 & 2006 & 29800 & $\sim$ & $\sim$ \Tstrut\Bstrut\\
		\hline
		\hline
	\end{tabular}
	\label{tab:run_time2}
\end{table}}

\subsection{Online Controller and Time-invariant Optimality Gap}\label{sec:result2}
To validate the optimality of proposed battery control policy in Theorem~\ref{theorem1} and the time-invariant gap in Theorem~\ref{theorem2}, we design nine test cases. Each test case has different market prices and battery round-trip efficiency $\eta=\eta\ud{d}\eta\ud{c}$ setting. In order to demonstrate the time-invariant property of the optimality gap, Case 7 to 9 are designed to double the duration of Case 4 to 6 with the same prices and efficiency setting. Each test case is simulated for 100 times using different randomly generated frequency regulation traces for reliability. At each time step, the signal is draw independently from a normal distribution with mean 0 and variance 1, and truncated between $[-1, 1]$. 

Table~\ref{tab:sim} summarizes the nine test case results. For each case, the penalty prices, round-trip efficiency, and the number of simulation control intervals used in each test case are listed, as well as the cycle depth bound $\hat{u}$ (calculated using \eqref{eq:pol3}) and the worst-case theoretical optimality gap $\epsilon$. Each time step is 1 minute. The maximum optimality gap incurs in 100 simulations are also recorded. 

This test validates Theorem~\ref{theorem2} since $\epsilon$ is exactly the same as the recorded maximum optimality gap of the proposed policy in all cases (both highlighted in pink). In particular, the proposed policy achieves exact control results in Case 1 to 3 because $\theta/\eta\ud{c}=\pi\eta\ud{d}$, while Case 4 to 9 have non-zero gaps because the round-trip efficiency is less than ideal ($\eta<1$). We also see that as penalty prices become higher, the optimal cycle depth $\hat{u}$ becomes larger and the battery follows the regulation instruction more accurately. Case 7 to 9 have the same parameter settings as to Case 4 to 6, except that the dimension of regulation signal doubled. The proposed policy achieves the same worst-case optimality gap in the two operation duration settings -- which again verifies that the worst-case optimality gap of the proposed online control policy is independent of operation time $T$.

\revise{\subsection{Comparison with Baseline Online Control Algorithms}
\begin{figure}[t]
	\centering
	\includegraphics[trim = 5mm 0mm 10mm 0mm, clip, width = .9\columnwidth]{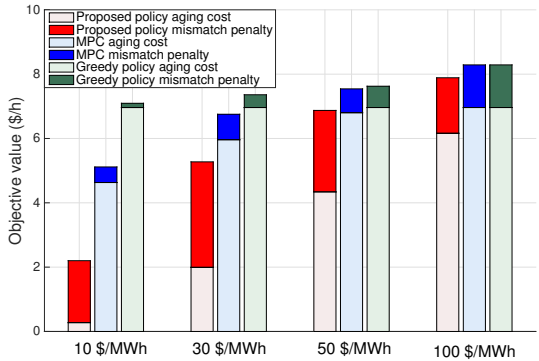}
	\caption{Regulation operating cost break-down comparison between the proposed policy, MPC and greedy policy. Although the proposed policy has higher mismatch penalties, the cost of cycle aging is significantly smaller, so it achieves better trade-offs between degradation and mismatch penalty.}
	\label{fig:pjm}
\end{figure}
\begin{figure}[t]
	\centering
	\includegraphics[trim = 5mm 0mm 10mm 0mm, clip, width = .9\columnwidth]{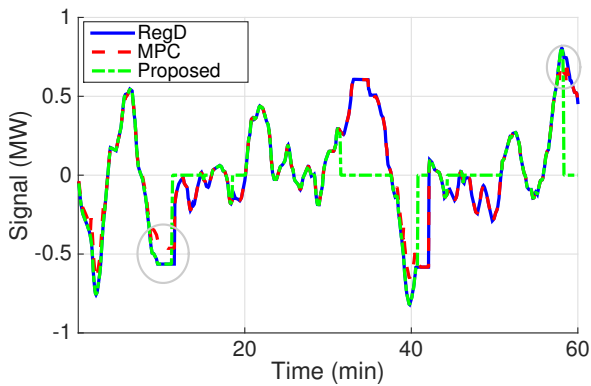}
	\caption{Comparison between the online MPC policy and proposed policy in response to the regulation signal RegD. Both methods have mismatch part (circled in gray) because of the tradeoff of degradation cost and mismatch penalty. MPC response is more aggressive because of under-estimate of the cycle depth.}
	\label{fig:pjm2}
\end{figure}
We take the proposed online battery control algorithm and compare against two state-of-the-art control algorithms for battery frequency regulation: the greedy control algorithm~\cite{bitar2011role} and the model predictive control (MPC) algorithm with rainflow degradation model~\cite{wang2017improving}. 

The greedy control strategy (in which BES follows the frequency regulation signal exactly within its power and energy limit) serves as the current market practice for BES frequency regulation [cite], and authors in~\cite{bitar2011role} showed it is optimal under linear energy throughput model. MPC is one of the most successful and the most popular real-time control methods. The basic idea of MPC is to predict the future instruction signal over a finite time horizon and compute the optimal control inputs. The obtained control is injected into the system until the next control step. Here, we assume that the future frequency regulation signal can be perfectly predicted within $60$ steps (4 minutes). Note, this is a strong assumption that does not hold in practice, since the second-by-second regulation signals are random and are almost impossible to predict~\cite{HughesEtAl2017}. Therefore, the following results indicate the best possible performance of such MPC-based methods.

We repeat the simulation using different penalty prices. We let $\theta=\pi$ in each test case and set the charging and discharging efficiency to $95\%$. Fig.~\ref{fig:pjm} summarizes the simulation results in the form of regulation operating cost versus penalty prices, the cycle aging cost and the regulation mismatch penalty are listed for each policy. Because the greedy control~\cite{bitar2011role} naively follow the signal and does not consider market prices, its control actions are the same in all price scenarios. The penalty increases linearly with the penalty price. MPC-based method tries to balance the degradation cost and mismatch penalty within a future look-ahead window. However, since the look-ahead window length is limited (due to real-time computational requirement and signal predictability), it tends to under-estimate the resulted cycle depth and leads to overaggressive control actions than the offline optima. Our proposed control policy leads to the lowest operating cost under all the cases by limiting the battery response within the optimal cycle depth. As the penalty price increases, the gap between the three policies becomes smaller since the optimal cycle depth $\hat{u}$ increases and battery follows the signal more closely. According to the historical billing data of PJM regulation market~\cite{PJM_regulation_RMPCP}, the regulation mismatch penalty price is usually below 50\$/MWh. Under such price setting, we save more than 30\% by using the proposed online control policy, and the battery can last as much as 3-4 times longer compared to the greedy and MPC controllers.}

\section{Conclusion and future work}
\label{sec:con}
We consider the optimal control of battery energy storage under a general ``pay-for-performance'' setup, where batteries need to trade-off between following instruction signals and the impact of degradation from charging and discharging actions.  We show that under electrochemically accurate cycle-based degradation models, the battery control problem can be formulated as a convex online optimization problem. Based on this result, we developed an online control policy that has a bounded time-invariant worst-case optimality gap, and is stricly optimal under certain market scenarios. From the case study in PJM regulation market, we verified the proposed degradation model and online control policy can significantly reduce operation cost and extend battery lifetime.

\revise{There are some natural directions for future work.
For example, the proposed threshold controller is the optimal under the pay for performance market settings, where the penalty prices are determined ahead of time. An important
future direction is to extend our results to settings where the
penalty prices are random in themselves, such as real-time
arbitrage~\cite{krishnamurthy20168, ParsonsEtAl2015}. Another interesting
direction to explore is online parameter estimation of battery degradation models. Currently, we use a fixed cycle stress function based on battery testing data from manufacturer. The performance of the proposed controller will be further enhanced if we could update the coefficients of the battery degradation model online, leveraging real-time battery condition measurements.}

\bibliographystyle{IEEEtran}
\bibliography{references,dc_bib}

\appendices
\section{Proof of Theorem \ref{theo1}}
\label{sec:appden1}
\setcounter{lemma}{0}
\setcounter{theorem}{0}
\setcounter{corollary}{0}
\subsection{Single step change convexity}
Here we continue proof of Theorem \ref{theo1} from Section \ref{sec:conv_initial}. Since all SoC profile can be written as the sum of step functions, by induction method, we first need to prove that $f(\bd x)$ is convex up to a step function as base case (Lemma \ref{lemma:single_step_proof}). Here, we first provide a restatement of the theorem and lemma that need to prove.
\revise{
\begin{theorem}\label{theorem1_appd}
	Suppose the battery cycle aging stress function $\Phi$ is convex. Then $f(\bd{x}) $, which is the mapping from the SoC profile $\bd{x}$ to the degradation cost:
	$$f(\bd{x}) = \left[\sum_{i=1}^{|\bd{v}|} \frac{\Phi(v_i)}{2} + \sum_{i=1}^{|\bd{w}|} \frac{\Phi(w_i)}{2}\right]\,,$$
	$$(\bd{v}, \bd{w}) = \text{Rainflow}(\bd{x})\,,$$
	is convex in terms of $\bd{x}$. That is, for any two SoC time series $\bd{x}, \bd{y} \in R^T$,
	$$f(\lambda \bd{x} + (1-\lambda) \bd{y}) \leq \lambda f(\bd{x}) + (1-\lambda) f(\bd{y}), \forall \lambda \in [0, 1]\,,$$
\end{theorem}
}
\revise{
	\begin{lemma} \label{lemma:single_step_proof2}
		Under the conditions in Theorem \ref{theo1}, the rainflow cycle-based cost function $f$ satisfies
		\begin{equation*}
		f\left(\lambda \mathbf{x} + (1-\lambda)Q_t U_t \right) \leq \lambda f(\mathbf{x}) + (1-\lambda) f(Q_t U_t )\,, \forall \; \lambda \in [0,1]\,,
		\end{equation*}
		where $\mathbf{x} \in \mathbb{R}^T$, and $Q_t U_t$ is a step function with a jump happens at time $t$ with amplitude $Q_t$.
\end{lemma}}

To show this, we need the following propositions.

%
\begin{proposition}\label{append_c1}
	Let $g(\cdot)$ be a convex function where $g(0) = 0$. Let $r_1$, $r_2$ be positive real numbers. Then
	\begin{equation}
	g(r_1+r_2) \geq g(r_1)+g(r_2)\,,
	\end{equation}
\end{proposition}

\begin{proof}
	By convexity of $g$, we have
	$$\frac{r_1}{r_1+r_2}g(r_1+r_2) + \frac{r_2}{r_1+r_2}g(0) \geq g(r_1)\,,$$
	and
	$$\frac{r_2}{r_1+r_2}g(r_1+r_2) + \frac{r_1}{r_1+r_2}g(0) \geq g(r_2)\,,$$
	Adding the two equations finish the proof.
\end{proof}


\begin{proposition}\label{append_c2}
	Let $g(\cdot)$ be a convex function where $g(0) = 0$. Let $r_1$, $r_2$ be positive real numbers, and $r_1 \geq r_2$. Then
	\begin{equation}
	g(r_1-r_2) \leq g(r_1)-g(r_2)\,,
	\end{equation}
\end{proposition}
\begin{proof}
	By Proposition \ref{append_c1}, $$g(\alpha+\beta) \geq g(\alpha)+g(\beta), \forall \alpha, \beta>0\,,$$
	Let $\alpha = r_1-r_2>0$, $\beta = r_2 >0$, so that
	$$g(r_1-r_2+r_2)  \geq g(r_1-r_2)+g(r_2).$$
\end{proof}
\begin{proposition}
	Let $g(\cdot)$ be a convex function where $g(0) = 0$. Let $r_1 \geq r_2 >0$ be positive real numbers. Then
	\begin{equation}
	g(\frac{1}{2} r_1- \frac{1}{2} r_2) \leq \frac{1}{2} g(r_1) - \frac{1}{2} g(r_2).
	\end{equation}
\end{proposition}
\begin{proof}
	From Proposition \ref{append_c2},
	$$g(\frac{1}{2} r_1- \frac{1}{2}r_2) \leq g(\frac{1}{2} r_1) - g(\frac{1}{2}r_2), \forall r_1 \geq r_2 >0$$
	Therefore, it suffices to show,
	$$g(\frac{1}{2} r_1) - g(\frac{1}{2} r_2) \leq \frac{1}{2} g(r_1) - \frac{1}{2} g(r_2)\,,$$
	Define $h(z) = g(\frac{1}{2} z)- \frac{1}{2} g(z)$,
	$$h'(z) = \frac{1}{2} g^{'}\big(\frac{1}{2} z\big) - \frac{1}{2} g\big(z\big) = \frac{1}{2} [g^{'}\big(\frac{1}{2} z\big)-g^{'}\big(z\big)]< 0\,,$$
	$h(\cdot)$ is a monotone decreasing function. For $r_1 \geq r_2>0$,
	\begin{align*}
	h(r_1) & \leq h(r_2)\,,\\
	g(\frac{1}{2} r_1)- \frac{1}{2} g(r_1) & \leq g(\frac{1}{2} r_2)- \frac{1}{2} g(r_2)\,,\\
	g(\frac{1}{2} r_1) - g(\frac{1}{2} r_2) & \leq \frac{1}{2} g(r_1) - \frac{1}{2} g(r_2).
	\end{align*}

	If $g(\cdot)$ is continous, we can generalize the midpoint property to a more broad $\lambda$,
	$$g(\lambda r_1- (1-\lambda) r_2) \leq \lambda g(r_1) - (1-\lambda) g(r_2), \forall \lambda r_1\geq  (1-\lambda) r_2 > 0$$

\end{proof}

\begin{proposition}
	Let $g(\cdot)$ be a convex function where $g(0) = 0$. Let $r_1, r_2, r_3$ be positive real numbers, which satisfy that $r_1+r_2-r_3 \geq 0$, and $r_i \leq r_1+r_2-r_3, \forall i \in [1,2,3]$. Then
	\begin{equation}
	g(r_1+r_2-r_3) \geq g(r_1)+g(r_2) - g(r_3)\,,
	\end{equation}
	\label{append_c4}
\end{proposition}
\begin{proof}
	From $r_1 \leq r_1+r_2-r_3$ we have $r_2 \geq r_3$. From $r_2 \leq r_1+r_2-r_3$, we have $r_1 \geq r_3$

	Let's further assume $r_1 \geq r_2$,
	$$g(r_1+r_2-r_3)-g(r_1) = (r_2-r_3) \cdot g^{'}(\theta_1), \theta_1 \in [r_1,r_1+r_2-r_3]\,,$$
	$$g(r_2)-g(r_3) = (r_2-r_3) \cdot g^{'}(\theta_2), \theta_2 \in [r_3,r_2]\,,$$
	Since $g(\cdot)$ is a convex function, for $\theta_2 \leq r_2 \leq r_1 \leq \theta_1$, we have $g^{'}(\theta_2) \leq g^{'}(\theta_1)$.
	Therefore,
	$$g(r_2)-g(r_3) \leq g(r_1+r_2-r_3)-g(r_1)\,,$$
	$$g(r_1+r_2-r_3) \geq g(r_1)+g(r_2)-g(r_3).$$

	If $r_1 < r_2$, similarly we have
	$$g(r_1)-g(r_3) \leq g(r_1+r_2-r_3)-g(r_2)\,,$$
	$$g(r_1+r_2-r_3) \geq g(r_1)+g(r_2)-g(r_3). $$

\end{proof}

\begin{proposition}
	Let $g(\cdot)$ be a convex function where $g(0) = 0$. Let $r_1, r_2, r_3, ..., r_n$ be real numbers, suppose
	\begin{itemize}
		\item $\sum_{i=1}^{n}r_i = D >0$
		\item $|r_i|  \leq D, \forall i \in \{1,2,3,...,n\}$
	\end{itemize}
	Then,
	\begin{equation}
	g(\sum_{i=1}^{n} r_i) \geq \sum_{\{i: r_i \geq 0\}} g(r_i)  - \sum_{\{i: r_i < 0\}} g(|r_i|).
	\end{equation}
	\label{append_c5}
\end{proposition}

\begin{proof}
	(1) If all $r_n$'s are positive, it is trivial to show $g(\sum_{i=1}^{n} r_i) \geq \sum_{i} g(r_i)$ by Proposition \ref{append_c1}.

	(2) If $r_n$ contains both positive and negative numbers, we order them in an ascending order and renumber them as,

	$$r_1 \leq r_2 \leq ... \leq 0 \leq ... \leq r_{n}\,,$$

	Pick $r_1$ (the most negative number), and find some positive $r_i$, $r_{i+1}$ such that
	$$r_{i+1} \geq r_{i} \geq |r_1| >0\,,$$

	Applying Proposition \ref{append_c4}, we have
	$$g(r_{i+1}+r_{i}+r_{1}) = g(r_{i+1} + r_{i} - |r_1|) \geq g(r_{i+1}) +g (r_{i}) - g(|r_1|)\,,$$

	Note, if we can not find such $r_{i}$, $r_{i+1}$, eg. $r_n < |r_1|$. We could group a bunch of \emph{postive} $r_i$'s to form two new variables $s_1 = \sum_{i \in N_1} r_i$, $s_2 = \sum_{i \in N_2} r_i$ where $N_1 \cap N_2 = \emptyset$.  For sure there exists such $s_1 \geq s_2 \geq |r_1|$, since
	\begin{align*}
	|\sum_{i=1}^{n} r_i| & = |\sum_{i: r_i \geq 0} r_i + \sum_{j: r_j < 0, j \neq 1} r_j + r_1| \\
	& =  \sum_{i: r_i \geq 0} r_i  - |\sum_{j: r_j < 0, j \neq 1} r_j| - |r_1|\\
	& = D
	\end{align*}
	$$\sum_{i: r_i \geq 0} r_i = |\sum_{j: r_j < 0, j \neq 1} r_j| + |r_1| + D \geq 2|r_1|$$

	Applying Proposition \ref{append_c4},
	\begin{align*}
	\ \  & \ g(s_{1}+s_{2}+r_{1}) \\
	= &\  g(s_{1} + s_{2} - |r_1|)\\
	\geq &\  g(s_{1}) +g (s_{2}) - g(|r_1|)\\
	\geq &\ \sum_{i \in N_1}  g(r_i) + \sum_{j \in N_2} g(r_j) - g(|r_1|), N_1 \cap N_2 = \emptyset
	\end{align*}

	Define $r_1^{'} = r_{i+1}+r_{i}+r_{1} > 0$ and re-order $r_1^{'}, r_2, r_3, ..., r_{i-1}, r_{i+2},...,r_{n}$. Or define ${r_1^{'}}= y_{1}+y_{2}+r_{1} > 0$, re-order $\big\{r_i: i \neq 1, i \notin N_1 \cup N_2\big\}, r_1^{'}$.

Repeat the above steps till all ${r_i}^{'}$ are postive and finish the proof.
\end{proof}

\begin{proposition}
	Consider a step change added to $\bd x$, where ${\bd x}^{'}(t) = \bd x(t)+Q_t U_t$, $t \in [0,T]$. Suppose $Q_t$ is positive \footnote{The proof for negative $Q_t$ is the same, just change $Q_t$ to $\left|Q_t\right|$}, the rainflow cycle decomposition results (only considering charging cycles) for $\bd x$ and ${\bd x}^{'}$ are,

	$$\bd x: v_1, v_2,..., v_m, ..., v_M\,,$$
	$${\bd x}^{'}: {v_1}^{'}, {v_2}^{'},..., {v_n}^{'}, ..., {v_N}^{'}\,,$$

	Define $L = max(M,N)$, we could re-write the cycles in $\bd x$ and ${\bd x}^{'}$ as,

	$$\bd x: \underbrace{v_1, v_2, ..., v_M, 0, 0,...}_{L}\,,$$
	$${\bd x}^{'}: \underbrace{{v_1}^{'}, {v_2}^{'}, ..., {v_N}^{'}, 0, 0,...}_{L}\,,$$

	Define $\Delta v_i$ such that,
	$${v_i}^{'} = v_i + \Delta v_i\,,\forall i = 1,2,...,L$$

	The following relations always holds,
	\begin{equation}
	\left|\sum_{i=1}^{L} \Delta v_i\right| \leq Q_t\,,
	\label{prop6:eq1}
	\end{equation}
	\begin{equation}
	|\Delta v_i| \leq Q_t\,,
	\label{prop6:eq2}
	\end{equation}
\label{append_c6}
\end{proposition}

\begin{proof}
There exists a small enough $\Delta Q$ such that only one cycle depth $v_i$ will change.

$$|\Delta v_i| \leq \Delta Q\,,$$
$$-\Delta Q \leq \Delta v_i \leq \Delta Q\,,$$

Consider $Q_i$ as a cumulation of small $\Delta Q$, by the principle of integration, we have
$$-\int \Delta Q dq \leq \sum_{i=1}^{L} \Delta v_i \leq \int \Delta Q dq\,,$$

Such that,
$$\left|\sum_{i=1}^{L} \Delta v_i\right| \leq Q_t$$

$|\Delta v_i| \leq Q_t$ holds for the worst case where all cycle depth changes happen at one certain cycle. Therefore, it is trivial to show that $|\Delta v_i| \leq Q_t$ hold in all conditions.

\end{proof}

By propositions 1-6, we get the proof of Lemma 1 below.

\begin{proof}
	Let's consider $\bd x^{'} = \lambda \bd x + (1-\lambda)Q_t U_t$. Then the rainflow cycle decomposition results for $\lambda \bd x$ and $\bd x^{'}$ are
	\begin{align*}
	\lambda \bd x: & \underbrace{\lambda v_1, \lambda v_2, ..., \lambda v_M, 0, 0,...}_{L} \\
	{\bd x}^{'}: & \underbrace{{v_1}^{'}, {v_2}^{'}, ..., {v_N}^{'}, 0, 0,...}_{L}
	\end{align*}
Define $\Delta v_i$ such that,
$${v_i}^{'} = \lambda v_i + (1-\lambda) \Delta v_i\,,\forall i = 1,2,...,L$$
\begin{small}
\begin{align}
	& f\big(\lambda \bd x + (1-\lambda)Q_t U_t\big) \nonumber\\
	= & \sum_{i=1}^{L} \Phi \big(\lambda v_i + (1-\lambda)\Delta v_i\big) \nonumber\\
	= & \sum_{i=1}^{L^{+}} \underbrace{\Phi \big(\lambda v_i + (1-\lambda)\Delta v_i\big)}_{\Delta v_i \geq 0}  + \sum_{i=1}^{L^{-}} \underbrace{\Phi \big(\lambda v_i - (1-\lambda)|\Delta v_i|\big)}_{\Delta v_i < 0} \nonumber\\
	\leq & \sum_{i=1}^{L^{+}} [\lambda \Phi(v_i) + (1-\lambda) \Phi(\Delta v_i)] + \sum_{i=1}^{L^{-}} [\lambda \Phi(v_i) - (1-\lambda) \Phi(|\Delta v_i|)] \nonumber\\
	\leq & \lambda \sum_{i=1}^{L} \Phi(v_i) + (1-\lambda) \big[\sum_{i=1}^{L^{+}} \Phi(\Delta v_i) - \sum_{i=1}^{L^{-}} \Phi(|\Delta v_i|)\big]
	\label{lemma1:eq1}
\end{align}
\end{small}

To continue the proof in \eqref{lemma1:eq1} and derive the final relation, we separate the whole variable space to two cases based on equations \eqref{prop6:eq1} and \eqref{prop6:eq2}.

\emph{(1). Assume $\sum_{i=1}^{L} \Delta v_i = Q_i$, $|\Delta v_i| \leq Q_i$}. By Proposition~\eqref{append_c5}, it follows that
\begin{align*}
& f\big(\lambda \bd x + (1-\lambda)Q_t U_t\big) \nonumber\\
\leq & \lambda \sum_{i=1}^{L} \Phi(v_i) + (1-\lambda) \big[\sum_{i=1}^{L^{+}} \Phi(\Delta v_i) - \sum_{i=1}^{L^{-}} \Phi(|\Delta v_i|)\big]\nonumber\\
\leq  & \lambda \sum_{i=1}^{L} \Phi(v_i) + (1-\lambda) \Phi(\sum_{i=1}^{L}  \Delta v_i)\nonumber\\
= & \lambda \sum_{i=1}^{L} \Phi(v_i) + (1-\lambda) \Phi(Q_t)
\end{align*}
\emph{(2) Assume $-Q_t \leq \sum_{i=1}^{L} \Delta v_i <Q_t$, $|\Delta v_i| \leq Q_t$}.

Add some ``virtual cycles'' $v_{L+1}^{'}, v_{L+2}^{'}, ..., v_{L+K}^{'}$ at the end of ${\bd x}^{'}$, each $v_{L+i}^{'}$ is positive and satisfies that $|v_{L+i}^{'}| \leq Q_t$. So that $\sum_{i=1}^{L+K} \Delta v_i = Q_t$, $|\Delta v_i| \leq Q_t, \forall i \in [1,2,...,L+K]$. Write 0 at the end of $\lambda \bd x$ to achieve the same cycle number.
\begin{align*}
	\lambda \bd x: & \underbrace{\lambda v_1, \lambda v_2, ..., \lambda v_M, 0, 0, 0,...,0}_{L+K} \\
	\bd x^{'}: & \underbrace{{v_1}^{'}, {v_2}^{'}, ..., {v_N}^{'}, 0, 0,...,0, v_{L+1}^{'}, v_{L+2}^{'}, ..., v_{L+K}^{'}}_{L+K}
\end{align*}

\resizebox{.97\linewidth}{!}{
	\begin{minipage}{\linewidth}
\begin{align*}
& f\big(\lambda \bd x + (1-\lambda)Q_t U_t\big) \nonumber\\
\leq &\lambda \sum_{i=1}^{L} \Phi(v_i) + (1-\lambda) \big[\sum_{i=1}^{l^{+}} \Phi(\Delta v_i) - \sum_{i=1}^{l^{-}} \Phi(|\Delta v_i|)\big]\nonumber\\
<& \lambda \sum_{i=1}^{L} \Phi(v_i) + (1-\lambda) \big[\sum_{i=1}^{l^{+}} \Phi(\Delta v_i) + \sum_{i=L+1}^{L+K} \Phi(\Delta v_i) - \sum_{i=1}^{l^{-}} \Phi(|\Delta v_i|)\big]\nonumber\\
\leq & \lambda \sum_{i=1}^{L} \Phi(v_i) + (1-\lambda) \Phi(\sum_{i=1}^{L+K}  \Delta v_i)\nonumber\\
= & \lambda \sum_{i=1}^{L} \Phi(v_i) + (1-\lambda) \Phi(Q_t)
\end{align*}
\end{minipage}
}

To sum up,

\begin{small}
\begin{align}
f\big(\lambda \bd x + (1-\lambda)Q_t U_t\big) & \leq \lambda \sum_{i=1}^{L} \Phi(v_i) + (1-\lambda) \Phi(Q_t) \nonumber\\
& = \lambda f(\bd x) + (1-\lambda) f(Q_t U_t)\,,
\end{align}
\end{small}
where $\lambda \in [0,1]$.
\end{proof}

Lemma 1 shows that $f(\bd x)$ is convex up to every step change in $\bd x$. Next, we will prove the general rainflow convexity by induction.

\subsection{General rainflow cycle life loss convexity}
We will prove the general rainflow convexity by induction.

By lemma 1, we already proved the base case convexity. When $K=1$,
$$f\big(\lambda \bd x + (1-\lambda)\bd y\big) \leq\lambda f(\bd x) + (1-\lambda) f(\bd y)\,, \lambda \in [0,1]$$

Next we need to show the induction relation. Suppose that, $f(\bd x)$ is convex up to the sum of $K$ step changes (arranged by time index)
$$f\big(\lambda \bd x + (1-\lambda) \bd y\big) \leq \lambda f(\bd x) + (1-\lambda) f(\bd y)\,, \lambda \in [0,1], \bd x, \bd y \in \mathbb{R}^{K}$$

\begin{figure}
	\centering
	\includegraphics[width=0.8 \columnwidth ]{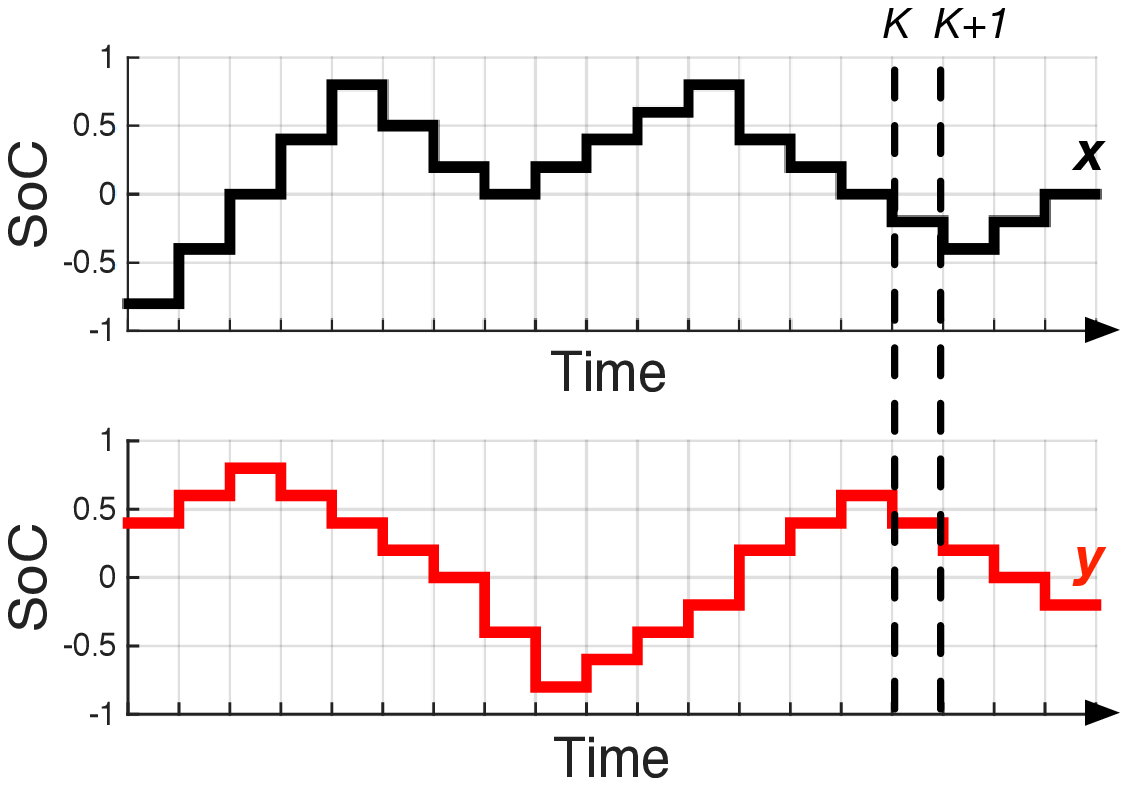}
	\caption{Induction step from $K$ to $K+1$, where we assume that $f$ is convex if one of the profiles (here the bottom profile) has only $K$ non-zero step changes and use it to show convexity of $f$ where one profile has $K+1$ non-zero steps.}
	\label{appd:fig_induction}
\end{figure}

Then we prove $f(\bd x)$ is convex up to the sum of $K+1$ step changes (see Fig. \ref{appd:fig_induction}),
$$f\big(\lambda \bd x + (1-\lambda) \bd y\big) \leq \lambda f(\bd x) + (1-\lambda) f(\bd y)\,, \lambda \in [0,1], \bd x, \bd y \in \mathbb{R}^{K+1}$$

The following proposition is needed for the proof.
\begin{proposition}\label{append_c7}
		\begin{equation}
		f(\sum_{t=1}^{K} P_t U_t) \geq f(\sum_{t=1}^{i-1} P_t U_t + (P_i+P_{i+1})U_i + \sum_{t=i+2}^{K} P_t U_t)\,,
		\end{equation}
	In other words, the cycle stress cost will reduce if combining adjacent unit changes.
\end{proposition}
\begin{proof}
	The rainflow cycle counting algorithm only considers local extreme points.

	I) If $P_i$ and $P_{i+1}$ are the same direction, combining them doesn't affect the value of local extreme points. Therefore the left side cost equals right side cost.

	II) If $P_i$ and $P_{i+1}$ are in different directions, suppose $P_i$ is negative and $P_{i+1}$ positive (otherwise the same). Time $t=i$ makes a local minimum point.
	\begin{itemize}
		\item Case a: If $|P_{i+1}| \leq |P_{i}|$, combining them will raise the value of local minimum point $i$, thus reducing the depth of cycles which contains $i$. Therefore, the cost after combining is less than the original cost.
		\item Case b: If $|P_{i+1}| > |P_{i}|$, combining them will lead to the removal of local minimum point $i$.

		In one case, if $P_{i-1}$ and $P_i$ are the same direction, time $t=i-1$ will make a local minimum point taking the place of time $t=i$. Therefore, the magnitude of the local minimum point decreases, similar to case (a), the total cost after combining is less than the original cost.

		In the other case, if $P_{i-1}$ and $P_i$ are different directions, we lose a full cycle with depth $|P_i|$ after combining. So the cost after combining $P_i$, $P_{i+1}$ is also less than the original.
	\end{itemize}
	To sum up, the cycle stress cost will reduce if combining adjacent unit changes.
\end{proof}
Recall the step function decomposition results for $\bd x$, $\bd y$ and $\lambda \bd x + (1-\lambda) \bd y$, where
\begin{align*}
\bd x = \sum_{t =1}^{T} P_t U_t, \bd y = \sum_{t =1}^{T} Q_t U_t \,, \\
\lambda \bd x + (1-\lambda) \bd y  = \sum_{t =1}^{T} Z_t U_t\,,
\end{align*}
There are three cases when $T$ goes from $K$ to $K+1$, classified by the value and symbols of $Z_{K}$, $Z_{K+1}$.

\noindent \textbf{Case 1:} $Z_{K}$ and $Z_{K+1}$ are same direction.

If $Z_{K+1}$ and $Z_{K}$ are same direction, we could move $Z_{K+1}$ to the previous step without affecting the total cost $f(\bd x + (1-\lambda) \bd y)$. Then we prove the $K+1$ convexity by applying Proposition \ref{append_c7}.
\begin{align}
	& f\left(\lambda \mathbf{x} + (1-\lambda) \mathbf{y}\right) \nonumber\\
	= & f\left(\lambda \mathbf{x}^{K} + (1-\lambda) \mathbf{y}^{K} + Z_{K+1}U_{K}\right) \nonumber\\
	= & f\left\{\lambda \mathbf{x}^{K} + (1-\lambda) \mathbf{y}^{K} + [\lambda P_{K+1} + (1-\lambda) Q_{K+1}] U_{K}\right\}\nonumber\\
	\leq & \lambda f\left(\mathbf{x}^{K}+ P_{K+1}U_{K}\right) + (1-\lambda)f\left(\mathbf{y}^{K}+ Q_{K+1}U_{K}\right) \nonumber\\
	\leq & \lambda f(\mathbf{x}) + (1-\lambda)f(\mathbf{y}) \ (\text{by Lemma \ref{lemma:single_step_proof}})
\end{align}
where $\mathbf{x}^{K}$ and $\mathbf{y}^{K}$ denote the first K elements of $\mathbf{x}$ and $\mathbf{y}$.

\noindent \textbf{Case 2:} $Z_{K}$ and $Z_{K+1}$ are different directions, with $|Z_{K}| \geq |Z_{K+1}|$. In this case, the last step $Z_{K+1}$ could be separated out from the previous SoC profile. Therefore,
\begin{small}
\begin{align}
	& f\left(\lambda \mathbf{x} + (1-\lambda) \mathbf{y}\right) & \nonumber\\
	= & f\left(\lambda \mathbf{x} + (1-\lambda) \mathbf{y}^{K}\right) + \Phi\left(Z_{K+1} U_{K+1}\right) \nonumber\\
	\leq & \lambda f(\mathbf{x}^{K}) + (1-\lambda)f(\mathbf{y}^{K}) + \Phi\left[\lambda P_{K+1} U_{K+1} + (1-\lambda) Q_{K+1}U_{K+1}\right] & \nonumber\\
	\leq & \lambda \big[f(\mathbf{x}^{K}) + \Phi(P_{K+1}U_{K+1})\big] + (1-\lambda)\big[f(\mathbf{y}^{K})  + \Phi(Q_{K+1} U_{K+1})\big] & \nonumber\\
	\leq &\lambda f(\mathbf{x}) + (1-\lambda)f(\mathbf{y}) &
\end{align}
\end{small}

\noindent \textbf{Case 3:} $Z_{K}$ and $Z_{K+1}$ are different directions, with $|Z_{K}| < |Z_{K+1}|$. In such condition, $Z_{K+1}$ is not easily separated out from previous SoC. To derive the induction relation, we analyze in three further sub-cases.
\begin{figure}
	\centering
	\includegraphics[width= 0.9 \columnwidth]{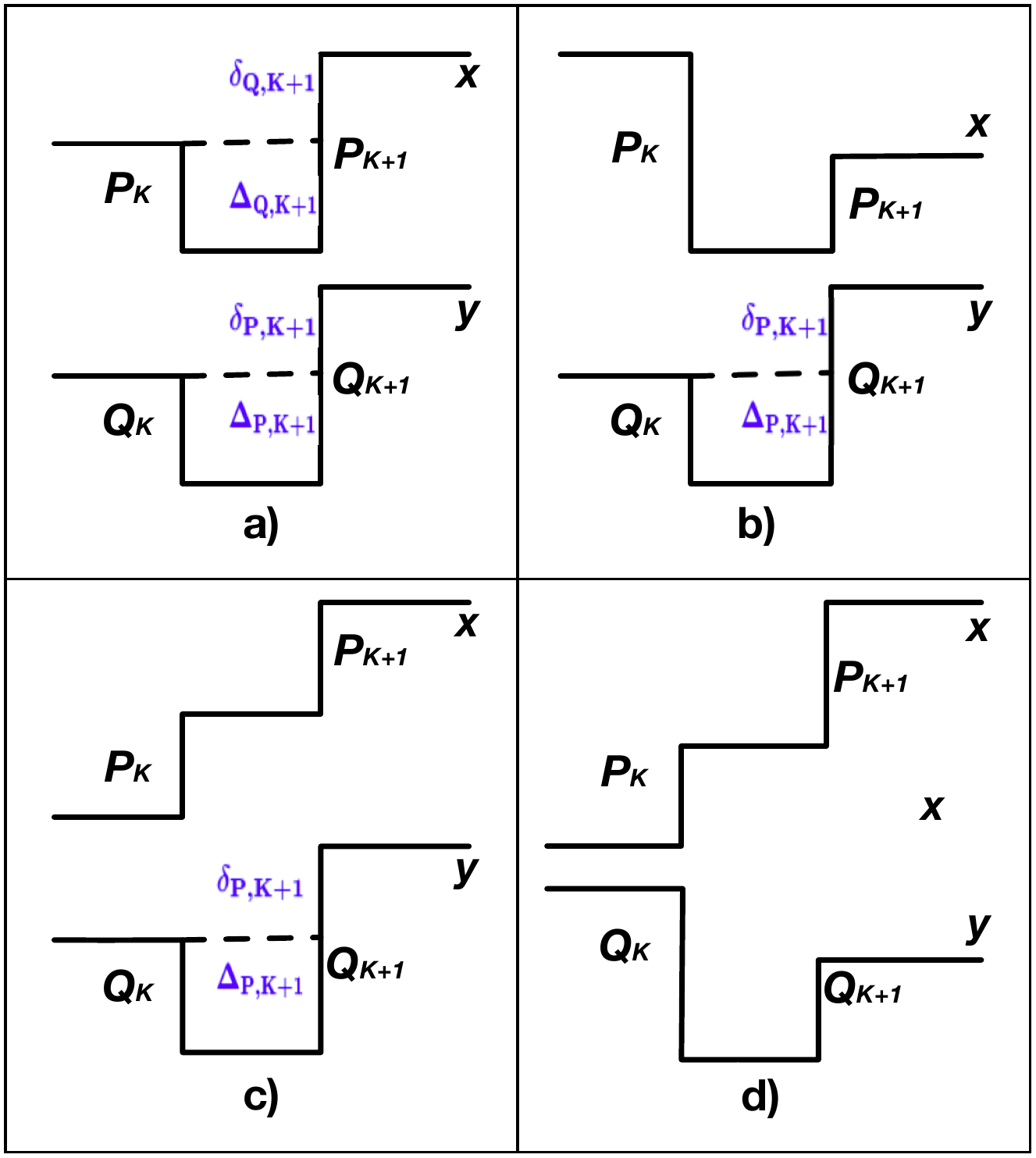}
	\caption{Four cases of $P_K$, $P_{K+1}$, $Q_K$, $Q_{K+1}$}
	\label{fig_append_4case}
\end{figure}
\begin{itemize}
	\item $Z_{K-1}$ and $Z_{K}$ are the same direction. In this sub-case, we could use the same ``trick" in Case 1 to combine step $K-1$ and $K$. Proof is trivial for this case.

	\item $Z_{K-1}$ and $Z_{K}$ are different directions, while $Z_{K}$ and $Z_{K+1}$ together form a cycle that is \emph{separable} from the rest of signal (eg. it is the deepest cycle). We can separate $Z_{K+1}$ out, and proof will be similar to Case 2.

	\item $Z_{K-1}$ and $Z_{K}$ are different directions, and  $Z_{K}$, $Z_{K+1}$ do not form a separate cycle. This condition is the most complicated case, since it's hard to move $Z_{K+1}$ to the previous step, or separate it out. Therefore, we need to look into $P_K$, $P_{K+1}$, $Q_K$, $Q_{K+1}$ in order to show the $K+1$ step convexity. It contains four more situations (Fig. \ref{fig_append_4case}), for simplicity we only consider the cost of charging cycles. Showing convexity for each situation finishes the overall convexity proof.

	Case a) $P_K$ and $P_{K+1}$ are in different directions, with $|P_{K}|<|P_{K+1}|$. $Q_K$ and $Q_{K+1}$ are also in different directions, with $|Q_{K}|<|Q_{K+1}|$. In such condition, the extra charging half cycle $\Delta_{K+1}$ could be decomposed into two charging half cycles in $\bd x$
	and $\bd y$ respectively.
	\begin{align}
	& f(\lambda \bd x + (1-\lambda) \bd y) \nonumber\\
 	= & f\left(\lambda \bd x^{K} + (1-\lambda) \bd y^{K} + Z_{K+1}U_{K+1}\right) \nonumber\\
 	= & f(\lambda \bd x^{K} + (1-\lambda)\bd y^{K} + Z_{K+1}U_{K}) + \Phi(\Delta_{K+1}) \nonumber\\
   	\leq & \lambda f\big(\bd x^{K}+P_{K+1}U_{K}\big) + (1-\lambda) f\big(\bd y^{K}+Q_{K+1}U_{K}\big) \nonumber\\
 	& + \Phi\big(\lambda \Delta_{P,K+1}+(1-\lambda) \Delta_{Q,K+1}\big) \nonumber\\
 	\leq & \lambda \big[f(\bd x^{K}+P_{K+1}U_{K})+\Phi(\Delta_{P,K+1})\big] \nonumber\\
 	& + (1-\lambda) \big[f(\bd y^{K}+Q_{K+1}U_{K})+\Phi(\Delta_{Q,K+1})\big] \nonumber\\
 	= & \lambda f(\bd x) + (1-\lambda) f(\bd y)
	\end{align}

	Case c) $P_K$ and $P_{K+1}$ are in the same direction. $Q_K$ and $Q_{K+1}$ are in different directions, with $|Q_{K}|<|Q_{K+1}|$.
	\begin{align}
 	& f(\lambda \bd x + (1-\lambda) \bd y) \nonumber\\
 	= & f(\lambda \bd x^{K} + (1-\lambda) \bd y^{K} + Z_{K+1}U_{K+1}) \nonumber\\
 	= & f(\lambda \bd x^{K} + (1-\lambda)\bd y^{K} + Z_{K+1}U_{K}) + \Phi(\Delta_{K+1}) \nonumber\\
 	\leq & \lambda f\big(\bd x^{K}\!+\!P_{K+1}U_{K}\big)\!+\! (1\!-\! \lambda) f\big(\bd y^{K}\!+\!Q_{K+1}U_{K}\big)\! +\! \Phi(\Delta_{K+1}) \nonumber\\
 	\leq & \lambda f\big(\bd x^{K}+P_{K+1}U_{K}\big) + (1-\lambda) f\big(\bd y^{K}+Q_{K+1}U_{K}\big) \nonumber\\
 	& + \Phi\big[(1-\lambda) \Delta_{P,K+1}-\lambda Q_{K+1}\big] \nonumber\\
 	\leq & \lambda f\big(\bd x^{K}+P_{K+1}U_{K}\big) + (1-\lambda) f\big(\bd y^{K}+Q_{K+1}U_{K}\big) \nonumber\\
 	& + \Phi\big[(1-\lambda) \Delta_{Q, K+1}\big] \nonumber\\
 	\leq & \lambda f(\bd x) + (1-\lambda) f(\bd y)
 	\end{align}

 	All the remaining task now is to prove the induction relation for cases b) and d).

 	Firstly, we note that b) implies d).  To show this, for case d), define $\hat{\bd x} = \sum_{t=1}^{K+1} \hat{P_t} U_t$ as a modified version of $\bd x$, where $\hat{P_t} = P_t$ for $t=1,...,K-1$, $\hat{P}_{K} = 0$, $\hat{P}_{K+1} = P_{K}+P_{K+1}$. We have $f(\hat{\bd x}) = f(\bd x)$. We also have that $f(\lambda \hat{\bd x}+(1-\lambda)\bd y) \geq f(\lambda \bd x+(1-\lambda)\bd y)$ because of the decreasing signal at $K$ for $\bd y$. Thus,
	    \begin{align}
	    f(\lambda \bd x + (1-\lambda)\bd y) &\leq f(\lambda \hat{\bd x}+(1-\lambda)\bd y) \nonumber\\
	   &  \stackrel{i)}\leq \lambda f(\hat{\bd x}) + (1-\lambda)f(\bd y) \nonumber\\
	   & = \lambda f(\bd x) +(1-\lambda)f(\bd y)
	    \end{align}

    where i) follows from assuming b) is true and letting $\Delta_{Q,K+1} =0$ and reversing the label of P and Q.

    Therefore we only need to prove case b). Case b) contains two different circumstances in terms of $Q_{K}$ and $Q_{K+1}$ contained in $\bd y$. We need the following proposition for the proof of case b).
    \begin{proposition}
	 Let $g$ be a convex increasing function and given real numbers $r_1 > r_2 >0$. Then $g(r_1)+g(r_2) \geq g(r_1-\delta) + g(r_2+\delta)$ if $r_2+\delta < r_1$, where $\delta$ is a small positive real number.
	\end{proposition}
    \begin{proof}
    	Define $h(x) = g(x)-g(x-\delta)$ and $x, x-\delta \geq 0$. We have,
    	$$h^{'}(x) = g^{'}(x)-g^{'}(x-\delta) \geq 0\,,$$
    	because $g$ is convex and $x > x-\delta$. Therefore, $h(\cdot)$ is an increasing function, $\forall a > b+\delta$,
    	\begin{align*}
    	h(r_1) & \geq h(r_2+\delta)\\
    	g(r_1) - g(r_1-\delta) &\geq g(r_2+\delta) - g(r_2)
    	\end{align*}
    	Moving $g(r_2)$ to the left side and $g(r_1-\delta)$ to the right side of the inequality finishes the proof.
    \end{proof}

	To not use too many negative signs, we denote $\bar{P}_t = -P_t$, $\bar{Q}_t = -Q_t$.

	\textcircled{1} $Q_{K}$, $Q_{K+1}$ do not form a cycle that is separate from the rest of $\bd y$.
	\allowdisplaybreaks
	\begin{align*}
	& f(\lambda \bd x + (1-\lambda)\bd y) \nonumber\\
	= & f(\lambda \bd x^{K} + (1-\lambda)\bd y^{K} + Z_{K+1}U_{K}) + \Phi(\Delta_{K+1}) \nonumber\\
	= & f(\lambda \bd x^{K} + (1-\lambda) \bd y^{K} + \lambda P_{K+1}U_{K} + (1-\lambda) Q_{K+1}U_{K})\nonumber\\
	& + f(\lambda \bar{P}_{K}U_{K+1} + (1-\lambda)  \bar{Q}_{K}U_{K+1}) \nonumber\\
	= & f\Big(\lambda \bd x^{K} + \lambda \bar{P}_{K}U_{K} - \lambda \bar{P}_{K}U_{K} + (1-\lambda) \bd y^{K} \nonumber\\
	& \ \ \ \ + \lambda P_{K+1} U_{K}  + (1-\lambda) Q_{K+1}U_{K}\Big) \nonumber\\
	& + f\Big(\lambda P_{K+1}U_{K+1} - \lambda P_{K+1}U_{K+1} + \lambda \bar{P}_{K}U_{K+1} \nonumber\\
	& \ \ \ \ + (1-\lambda) \bar{Q}_{K}U_{K+1}\Big) \nonumber\\
	= & f\Big(\lambda \bd x^{K} + \lambda \bar{P}_{K}U_{K}  +(1-\lambda)\big(\bd y^{K}+ Q_{K+1}U_{K} \nonumber\\
	&\ \ \ \ + \frac{\lambda}{1-\lambda}(P_{K+1}U_{K}- \bar{P}_{K}U_{K}) \big)\Big) \nonumber\\
	&+ f \Big(\lambda P_{K+1}U_{K+1} + (1-\lambda) \big( \bar{Q}_{K}U_{K+1}\nonumber\\
	&\ \ \ \ +\frac{\lambda}{1-\lambda}(\bar{P}_{K}U_{K+1}- P_{K+1}U_{K+1} )\big)\Big) \nonumber\\
	\leq & \lambda f(\bd x^{K}+\bar{P}_{K}U_{K}) + \lambda f(P_{K+1}U_{K+1} ) \nonumber\\
	&+(1-\lambda)f\big(\bd y^{K}+\big(Q_{K+1}-\frac{\lambda}{1-\lambda}(\bar{P}_{K}-P_{K+1})\big)U_{K}\big)\nonumber\\
	&+(1-\lambda)f\big(\big(\bar{Q}_{K}+\frac{\lambda}{1-\lambda}(\bar{P}_{K}-P_{K+1})\big)U_{K+1}\big)
\end{align*}

Only considering charging cycles, the first line is the cost of $\lambda f(\bd x)$. Now we show,
\begin{align*}
f\big(\bd y^{K}+\big(Q_{K+1}-\frac{\lambda}{1-\lambda}(\bar{P}_{K}-P_{K+1})\big)U_{K}\big)\nonumber\\
+f\big(\big(\bar{Q}_{K}+\frac{\lambda}{1-\lambda}(\bar{P}_{K}-P_{K+1})\big)U_{K+1}\big) \leq f(\bd y)\,,
\end{align*}

We can write out the cost of charging cycles in $\bd y$ as $f(\bd y)$,
\begin{equation*}
f(\bd y) = \sum_{i}^{N-1} \Phi(v_i) + \Phi(\bar{Q}_{K})\,,
\end{equation*}

where $v_{N-1}$ is the charging cycle that the $K+1$ step of $\bd y$, denoted as ${\bd y}(K+1)$ belongs to.

By assumption that $Q_{K}$ and $Q_{K+1}$ do not form a separate cycle, so that $v_{N-1} \geq Q_{K+1}$. By assumption, $\frac{\lambda}{1-\lambda}(\bar{P}_{K}-P_{K+1}) > 0$. Let  $\delta=\frac{\lambda}{1-\lambda}(\bar{P}_{K}-P_{K+1})$, and since $P_{K+1}+Q_{K+1} \geq \bar{P}_{K}+\bar{Q}_{K}$ in case b), $\bar{Q}_{K}+\delta \leq Q_{K+1} \leq v_{N-1}$. Therefore applying Proposition 8, $a = v_{N-1}$ and $b = \bar{Q}_{K}$, $\delta=\frac{\lambda}{1-\lambda}(\bar{P}_{K}-P_{K+1})$, we have the desired result.

\textcircled{2} The other case is that $Q_{K}$, $Q_{K+1}$ form a cycle that is separate from the rest of $\bd y$. Similar as case \textcircled{1}, we need to show
\begin{align*}
f\big(\bd y^{K}+\big(Q_{K+1}-\frac{\lambda}{1-\lambda}(\bar{P}_{K}-P_{K+1})\big)U_{K}\big) & \nonumber\\
+f\big(\big(\bar{Q}_{K}+\frac{\lambda}{1-\lambda}(\bar{P}_{K}-P_{K+1})\big)U_{K+1}\big)\leq f(\bd y)&
\end{align*}

Since $Q_{K+1} = \bar{Q}_{K}+\delta_{Q, K+1}$, we re-write the above inequality as,
\begin{align*}
f\big(\bd y^{K}+\big(\bar{Q}_{K}+\delta_{Q, K+1}-\frac{\lambda}{1-\lambda}(\bar{P}_{K}-P_{K+1})\big)U_{K}\big)& \nonumber\\
\!+\!f\big(\big(Q_{K\!+\!1}\!-\!\delta_{Q,\! K+1}\!+\! \frac{\lambda}{1\!-\!\lambda}(\bar{P}_{K}\!-\!P_{K+1})\big)U_{K\!+\!1}\big) \leq\!f\!(\bd y\!)&
\end{align*}

Denote $\delta = \delta_{Q, K+1}-\frac{\lambda}{1-\lambda}(\bar{P}_{K}-P_{K+1}) > 0$
\begin{equation*}
f(\bd y) = \sum_{i}^{N-2} \Phi(v_i) + \Phi (v_{N-1})+ \Phi(Q_{K+1})\,,
\end{equation*}

$v_{N-1}$ is the deepest charging cycle where its ending SoC equals to $Q_{K}$'s starting SoC, and $v_{N-1} \leq Q_{K+1}$ since $Q_{K+1}$ forms a separate cycle. Applying Proposition 8 by setting $a = Q_{K+1}$, $b = v_{N-1}$, we have the desired result.
\end{itemize}

\section{Proof of Theorem \ref{theorem2}}
\label{sec:appden2}
\revise{
\begin{theorem}[ \textbf{Online optimality}] \label{opt_theorem2}
	\noindent \emph{Suppose the battery cycle aging stress function $\Phi$ is strictly convex. The proposed online control policy $g$
		in Algorithm~\ref{alg:online} Section~\ref{sec:policy}. A, has a constant worst-case optimality gap that is independent of the operation time duration $T$.}
	\begin{align}
		\sup_{\bd r} (J_{g}-J^*) \leq \epsilon, \text{$\forall$ $x_0$ and  $\forall$ $\{r_t\}$, $t\in\{1,\dots,T\}$.} \nonumber
	\end{align}
\end{theorem}
\begin{corollary}[\textbf{Zero-optimality Gap}]
	\noindent\emph{If $\pi\eta\ud{d} = \theta/\eta\ud{c}$, then there is no gap between the proposed online algorithm and the optimal value of solving the offline problem (given entire $\bd r \in \mathbb{R}^T$).}
	\label{opt:gap_0}
\end{corollary}
\begin{theorem}\label{theorem3}
	If function $\Phi(\cdot)$ is strictly convex, then the worst-case optimality gap for the proposed policy $g(\cdot)$ in Theorem \ref{theorem1} is
	\begin{align}
	\epsilon = \begin{cases}
	\epsilon\ud{w} & \text{if $\pi\eta\ud{d} > \theta/\eta\ud{c}$} \\
	0 & \text{if $\pi\eta\ud{d} = \theta/\eta\ud{c}$} \\
	\epsilon\ud{v} & \text{if $\pi\eta\ud{d} < \theta/\eta\ud{c}$}
	\end{cases}\,
	\end{align}
	where
	\begin{align}
	\epsilon\ud{w} &= J\ud{w}(\hat{u})+2J\ud{v}(\hat{u})-J\ud{w}(\hat{w})-2J\ud{v}(\hat{v}) \label{eq:gap_v2}\\
	\epsilon\ud{v} &= 2J\ud{w}(\hat{u})+J\ud{v}(\hat{u})-2J\ud{w}(\hat{w})-J\ud{v}(\hat{v})\label{eq:gap_u2}\,.
	\end{align}
\end{theorem}}
\subsection{Model Reformulation}
Both Theorem \ref{theorem1} and Corollary \ref{thm:gap_0} follow directly from Theorem~\ref{theorem2}. To prove Theorem~\ref{theorem2}, we rewrite the optimization problem \eqref{eq:battery_operation2} in Section III.E as,
\begin{subequations}\label{eq:opt2}
\begin{align}
    (\mathbf{c}^*, \mathbf{d}^*) \in \mathrm{arg} &\min_{\mathbf{c}, \mathbf{d}} f(\mathbf{c}, \mathbf{d}) - \tau \sum_{t=1}^T\big[\theta c_t + \pi d_t\big]\label{eq:the1_0}\\
    &\text{subject to \eqref{Eq:PF_C4}, \eqref{Eq:PF_C5}, and} \nonumber\\
    & 0\leq c_t \leq [r_t]^+  \label{eq:the1_1}\\
    & 0\leq d_t \leq [-r_t]^+ \label{eq:the1_2}
\end{align}
\end{subequations}
by observing that a battery's actions would never exceed the regulation signals. $f(\mathbf{c}, \mathbf{d})$ defines the rainflow cycle-based degradation cost.

%

We utilize the rainflow algorithm to transform the problem into a cycle-based form. The rainflow method maps the entire operation uniquely to cycles, the sum of all charge and discharge power can be represented as the sum of cycle depths as (recall that a full cycle has symmetric depth for charge and discharge)
\begin{align}
    \sum_{i=1}^{|\mathbf{u}|}u_i + \sum_{i=1}^{|\mathbf{v}|}v_i &= \frac{\tau \eta\ud{c}}{E}\sum_{t=1}^T c_t
    \label{eq:rf1}\\
    \sum_{i=1}^{|\mathbf{u}|}u_j + \sum_{i=1}^{|\mathbf{w}|}w_i &= \frac{\tau}{\eta\ud{d}E}\sum_{t=1}^T d_t\,.\label{eq:rf2}
\end{align}
We substitute \eqref{eq:rf1} and \eqref{eq:rf2} into the reformulated objective function \eqref{eq:the1_0} to replace $c_t$ and $d_t$ with cycle depths
\begin{align}
    &J\ud{cyc}(\mathbf{c}, \mathbf{d}) + J\ud{reg}(\mathbf{c}, \mathbf{d}, \mathbf{r}) =\nonumber\\ &\sum_{i=1}^{|\mathbf{u}|}J\ud{u}(u_i) + \sum_{i=1}^{|\mathbf{v}|}J\ud{v}(v_i) + \sum_{i=1}^{|\mathbf{w}|}J\ud{w}(w_i)\,.
\end{align}

\subsection{Proof of Theorem~\ref{theorem2}}

The following lemmas support the proof of Theorem~\ref{theorem2}.
\begin{lemma}\label{lemma1}
	Suppose an minimizer $(\mathbf{c}^*, \mathbf{d}^*)$ of \eqref{eq:battery_operation2} in the offline setting has the corresponding cycle depths $(\mathbf{u}^*, \mathbf{v}^*, \mathbf{w}^*)$. Then the depth of each cycle in this result either reaches the optimal cycle depth or bounded by the operation constraints as
	\begin{subequations}
		\begin{align}
		u_i^* &= \min(\hat{u}, \overline{u}_i) \\
		v_i^* &= \min(\hat{v}, \overline{v}_i) \\
		w_i^* &= \min(\hat{w}, \overline{w}_i)
		\end{align}
	\end{subequations}
	where $\overline{u}_i$, $\overline{v}_i$, $\overline{w}_i$ denote constraint bounds including the regulation instruction signal and battery energy limit.
\end{lemma}

\begin{lemma}\label{lemma2}
	A cycle depth in the control action of $g(\cdot)$ either reaches the depth of $\hat{u}$ or is bounded by the operation constraints.
\end{lemma}

\begin{lemma}\label{lemma3}
	There exists one and only one half cycle with the largest depth in a rainflow residue profile. Other half cycles are in strictly decreasing order either to the left- or to the right-hand side direction of this largest half cycle.
	\begin{figure}[ht]
		\centering
		\includegraphics[trim = 5mm 0mm 10mm 0mm, clip, width = .95\columnwidth]{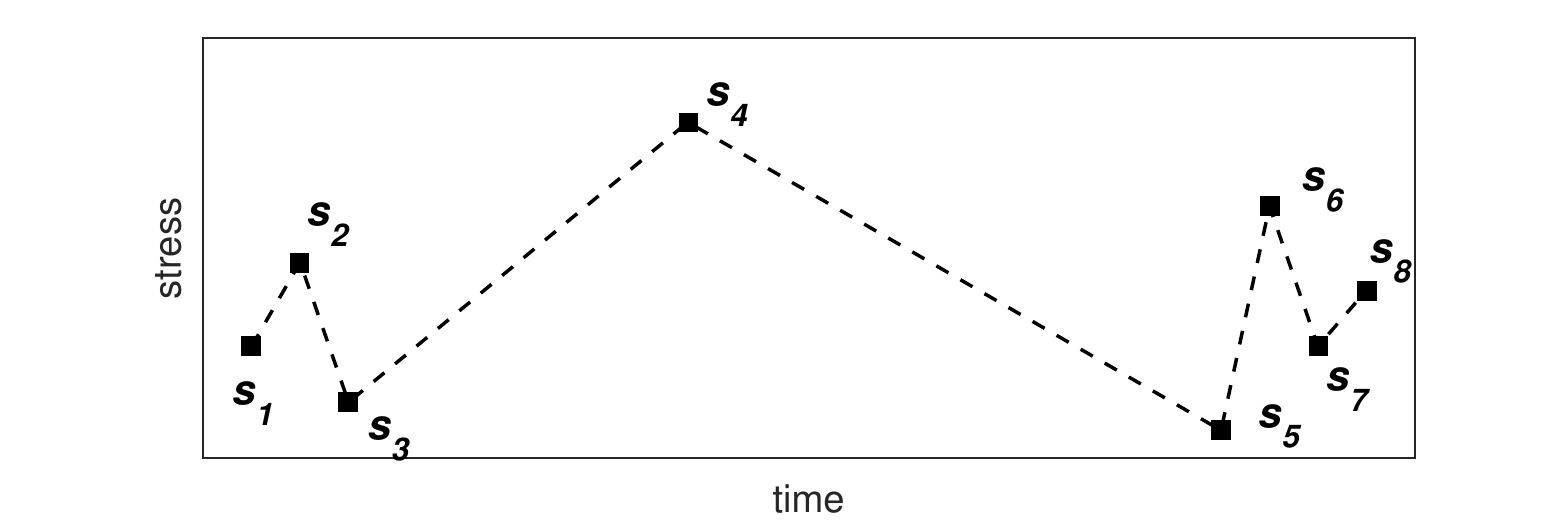}
		\caption{Illustration for Lemma~\ref{lemma3}. The largest half cycle is between $s_4$ and $s_5$, other half cycles are in strictly decreasing order either to the left- or to the right-hand side direction of this largest half cycle.}
		\label{fig:lemma3}
	\end{figure}
\end{lemma}

It is easy to see now from Lemma~\ref{lemma1} and Lemma~\ref{lemma2} that the proposed control policy achieves optimal control result for all full cycles, and the optimality gap is caused by half cycle results. Consider the following relationship in a rainflow residue profile as in Lemma~\ref{lemma3} assuming the largest half cycle is in the discharging direction
\begin{align}\label{the3:eq2}
\dotsc < w^*_{j-1} < v^*_{j-1} <w^*_j>v^*_j >w^*_{j+1}>\dotsc
\end{align}
and substitute Lemma~\ref{lemma2} into \eqref{the3:eq2}
\begin{align}
\dotsc \min\{\hat{v}, \overline{v}_j\} <\min\{\hat{w}, \overline{w}_j\}>\min\{\hat{v}, \overline{v}_{j+1}\}\dotsc
\end{align}
It is easy to see now that if $\hat{w}>\hat{v}$, then the largest possible value for $ w^*_j$ is $\hat{w}$,
and the largest possible value for $v^*_j$ and $v^*_{j-1}$ is $\hat{v}$, the rest half cycles in \eqref{the3:eq2} must have depths smaller than $\hat{v}$, which indicates that their depths are bounded by operation. If $\hat{v}>\hat{w}$, then the largest possible value for $ w^*_j$ is $\hat{w}$, and the rest half cycles must have depths smaller than $\hat{w}$. We repeat this analysis for cases that $v^*_j$ is the largest cycle, and summarize the half cycle conditions in Table~\ref{tab:half_cycle}
\begin{table}[ht]
	\centering
	\caption{Summarizing Half Cycle Depth Conditions}
	\begin{tabular}{l c c}
		\hline
		\hline
		&  $\hat{w}>\hat{v}$ & $\hat{w}<\hat{v}$ \Tstrut\Bstrut\\
		\hline
		Half cycles of depth $\hat{w}$ & At most one & At most two\Tstrut\Bstrut\\
		Half cycles of depth $\hat{v}$ & At most two & At most one\Tstrut\Bstrut\\
		Rest half cycles & must be $<\hat{v}$ & must be $<\hat{w}$\Tstrut\Bstrut\\
		\hline
		\hline
	\end{tabular}
	\label{tab:half_cycle}
\end{table}
Hence, the worst-case optimality gap is caused by that some half cycles have depth $\hat{u}$ or $\hat{w}$, while the control policy enforces $\hat{u}$ as the depth of all cycles unbounded by operation. The gap in Theorem~\ref{theorem2} is therefore calculated using half cycle depth conditions in Table~\ref{tab:half_cycle}.

\noindent \emph{Proof of Lemma~\ref{lemma1}:}
Since cycles are linear combinations of charge and discharge power, and constraints \eqref{eq:the1_1}, \eqref{eq:the1_2}, \eqref{Eq:PF_C5} can be transformed into linear constrains with respect to cycle depths. From Theorem~\ref{theo1}, the transformed cycle-based problem is also has a convex objective function with linear constraints. Although exact formulations of the transformed constraints are complicated to express, we use $\overline{u}_i$, $\overline{v}_i$, and $\overline{w}_i$ to denote these binds, which are sufficient for the proof of Theorem~\ref{theorem2}.\\

\noindent \emph{Proof of Lemma~\ref{lemma2}:}
\begin{figure}[ht]
    \centering
    \includegraphics[trim = 5mm 0mm 10mm 0mm, clip, width = .95\columnwidth]{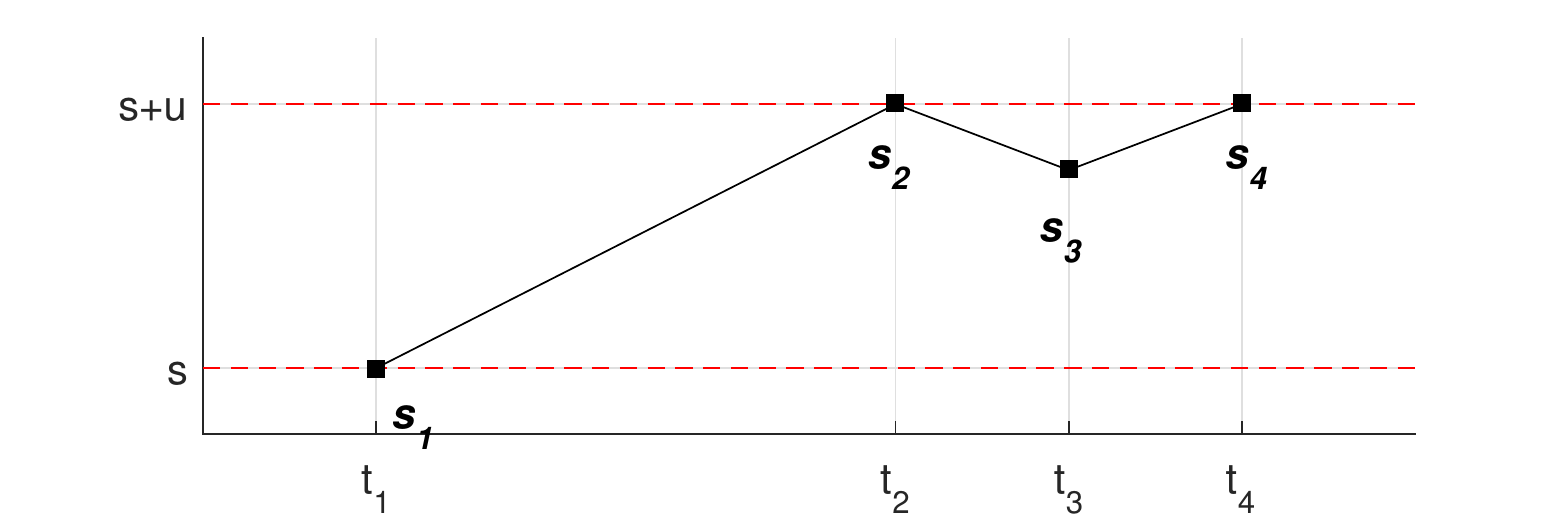}
    \caption{Illustration for Lemma~\ref{lemma2}.}
    \label{fig:lemma2}
\end{figure}
The rainflow method always identify the largest cycle as between the minimum and the maximum SoC point. In the proposed policy, any operation that goes outside the defined operation zone will cause the largest cycle depth to change instead of the depth of the cycle it was previous in. For example, in Fig.~\ref{fig:lemma2} the maximum cycle is between SoC $s$ and $s+u$, and the battery is at time $t_4$. If the battery continue to charge and the SoC goes about $s+u$, then this operation will increase the largest cycle depth instead of the shallower cycles assoicated with extrema $s_2$, $s_3$ and $s_4$.\\

\noindent \emph{Proof of Lemma~\ref{lemma3}:}
Because the rainflow method identifies a cycle from extrema distances if $\Delta s_{i-1}\geq \Delta s_{i} \leq \Delta s_{i+1}$, then all extrema in the rainflow residue must satisfy either $\Delta s_{i-1}< \Delta s_{i}  < s_{i+1}$ or $\Delta s_{i-1}< \Delta s_{i}  > s_{i+1}$ or $\Delta s_{i-1}> \Delta s_{i} > s_{i+1}$, which proofs this lemma.



\end{document}